\documentclass[a4paper,12pt,oneside,reqno]{amsart}

\usepackage[latin1]{inputenc}
\usepackage{hyperref}
\usepackage[headinclude,DIV13]{typearea}
\areaset{16cm}{24.7cm}
\usepackage{txfonts,amssymb,amsmath,amsthm,bbm,wasysym}
\usepackage[pdftex]{graphicx}
\usepackage{geometry}
\usepackage[framemethod=tikz]{mdframed}
\usepackage{cancel} 
\usepackage{graphicx} 
\graphicspath{ {figures/} }
\usepackage{caption}
\usepackage{subcaption}
\usepackage{mathrsfs}
\usepackage{xcolor}
\usepackage{enumerate}
\usepackage{nicefrac}	
\theoremstyle{plain}
\newtheorem{theorem}{Theorem}[section]
\newtheorem{lemma}[theorem]{Lemma}
\newtheorem{proposition}[theorem]{Proposition}
\newtheorem{corollary}[theorem]{Corollary}

\newtheorem{example}{Example}

\theoremstyle{definition}
\newtheorem{definition}[theorem]{Definition}
\newtheorem{assumption}{Assumption}

{\itshape}{\rmfamily}
{\itshape}{\rmfamily}

\theoremstyle{remark}
\newtheorem{remark}{Remark}

\def\paragraph#1{\noindent \textbf{#1}}

\numberwithin{equation}{section}

\def\d{\mathrm{d}}

\def\<{\langle}
\def\>{\rangle}

\def\e{\e}

\def\f{\phi}

\def\s{\sigma}

\def\R{{\Bbb R}}  
\def\N{{\Bbb N}}  %
\def\P{{\Bbb P}}  
\def\Z{{\Bbb Z}}  

\let\cal=\mathcal

\def\GG{{\cal G}}

\def\LL{{\cal L}}
\def\MM{{\cal M}}

 \def \e {{\epsilon}}
 \def \s {{\sigma}}

 \def \d {{\delta}}

 \def \f {{\phi}}

 \newcommand{\be}{\begin{equation}}
 \newcommand{\ee}{\end{equation}}

\newcommand{\bea}{\begin{eqnarray}}
 \newcommand{\eea}{\end{eqnarray}}
 
\def\TH(#1){\label{#1}}\def\thv(#1){\ref{#1}}
\def\Eq(#1){\label{#1}}\def\eqv(#1){(\ref{#1})}

 \def \1{\mathbbm{1}}

\def\eps{\varepsilon}
\def\v{\textbf{v}}

\def\ee{\mathrm{e}}
\DeclareMathOperator*{\argmin}{arg\,min}

\begin{document}
 \title[]{Stochastic individual-based models with power law mutation rate on a general finite trait space}
  
\author[]{Loren Coquille, Anna Kraut, Charline Smadi}
\address{L.\ Coquille\\Univ. Grenoble Alpes, CNRS, Institut Fourier, 38610 Gi\`eres, France}
\email{Loren.Coquille@univ-grenoble-alpes.fr}
\address{A.\ Kraut\\ Institut f\"ur Angewandte Mathematik\\
Rheinische Friedrich-Wilhelms-Universit\"at\\ Endenicher Allee 60\\ 53115 Bonn, Germany}
\email{kraut@iam.uni-bonn.de}
\address{C.\ Smadi\\ Univ. Grenoble Alpes, INRAE, LESSEM, F-38402 St-Martin-d'H\`eres, France
 and Univ. Grenoble Alpes, CNRS, Institut Fourier, 38610 Gi\`eres, France}
\email{charline.smadi@inrae.fr}



\thanks{This work was partially supported by the Deutsche Forschungsgemeinschaft (DFG, German Research Foundation) under Germany's Excellence Strategy GZ 2047/1, Projekt-ID 390685813 and GZ 2151, Project-ID 390873048 and through the Priority Programme 1590 ``Probabilistic Structures in Evolution''.
This work was also partially funded by the Chair \emph{Mod\'elisation Math\'ematique et Biodiversit\'e} of VEOLIA-Ecole Polytechnique-MNHN-F.X.
The authors thank Anton Bovier for stimulating discussions at the beginning of this work and comments, and L. Coquille and C. Smadi also thank him for his invitations
and welcome at Bonn University. The authors would also like to thank Sylvain Billiard for pointing out some bibliographical references.}

\begin{abstract}
	We consider a stochastic individual-based model for the evolution of a haploid, asexually reproducing population. 
	The space of possible traits is given by the vertices of a (possibly directed) finite 
	graph $\GG=(V,E)$. 
	The evolution of the population is driven by births, deaths, competition, and mutations along the edges of $\GG$. 
	We are interested in the large population limit under a mutation rate $\mu_K$ given by a negative power of the carrying capacity $K$ of the 
	system: $\mu_K=K^{-1/\alpha},\alpha>0$.
	This results in several mutant traits being present at the same time and competing for invading the resident population.
	We describe the time evolution of the orders of magnitude of each sub-population on the $\log K$ time scale, as $K$ tends to infinity. 
	Using techniques developed in \cite{htrans19}, we show that these are piecewise affine continuous functions, whose slopes are given by an algorithm 
	describing the changes in the fitness landscape due to the succession of new resident or emergent types. This work generalises 
	\cite{BovKra2018} to the stochastic setting, and Theorem 3.2 of \cite{BovCoqSma2018} to any finite mutation graph.
	We illustrate our theorem by a series of examples describing surprising phenomena arising from the geometry of the graph and/or the rate of mutations.
\end{abstract}

\maketitle

\section{Introduction}
\emph{Adaptive dynamics} is a biological theory that was developed to study the interplay between ecology and evolution. It involves the three mechanisms of heredity, 
mutations, and natural selection. It was first introduced in the 1990ies by Metz, Geritz, Bolker, Pacala, Dieckmann, Law, and coauthors 
\cite {MG96,DL96,GMK1997,BolPac1,BolPac2,DieLaw}, who mostly considered a deterministic setting but also heuristically mentioned first stochastic versions.
A paradigm of adaptive dynamics is the separation of the slow evolutionary and the fast ecological time scales, which is a result of reproduction with rare mutations. 
Invasion, fixation or extinction of a mutant population is determined by its invasion fitness, that describes the exponential growth rate of a single mutant 
in the current (coexisting) population(s) at equilibrium.

Stochastic individual-based models of adaptive dynamics have been 
rigorously constructed and first studied in the seminal work of Fournier and M\'{e}l\'{e}ard \cite{FM04}, and there is now a growing literature on these models. 
The population consists of a collection of individuals who reproduce, with or without mutation, or die after random exponential times depending 
on the current state of the whole population. The population size is controlled by a \emph{carrying capacity} $K$ which represents the amount of available 
resources.  
This class of models has first been studied in the original context of separation between evolutionary and ecological time scales. That is in the joint limit of large populations and rare mutations such that a mutant either dies out or fixates before the next mutation occurs. Mathematically this amounts to considering a probability of mutation satisfying in particular
\begin{equation}\label{mut-champagnat}
	\mu_K\ll1/{K\log K} \quad\text{ as } K\to\infty.
\end{equation}
We will call this regime 'rare mutation regime' in the sequel.
The description of the succession of mutant invasions, on the mutation time scale $1/K\mu_K$, in a monomorphic \cite{Cha06} or polymorphic \cite{CM11,BilSma15} asexual population gives rise respectively to the so-called \emph{Trait Substitution Sequence} or \emph{Polymorphic Evolution Sequence}. 
Extensions of the question to sexual populations were then studied, both in the haploid \cite{smadi2015eco,Coron2017} and the diploid \cite{CMM13,BovNeu} cases.

It is natural to consider the effect of a higher mutation rates, where mutation events are no longer separated, if we want to describe several mutant traits being present microscopically at the same time and competing for invading the resident population. The mutation rate given by 
\begin{equation}\label{mutation-rate}
\mu_K=K^{-\frac{1}{\alpha}},\quad\text{ for } \alpha>0
\end{equation} 
was considered in different contexts \cite{DurMay2011,Smadi2017,BovCoqSma2018,htrans19} 
and will be the concern of the present paper.
Notice that another mutation scale has been considered in \cite{BilSma15,BilSma2019} to model the interaction of few mutants in the case without recurrent mutations, namely 
$\mu_K$ of order $1/K\log K$.
\\

Another approach to adaptive dynamics has been introduced by Maynard Smith \cite{MS1970} under the name of \emph{adaptive walks}. This was further developed by 
Kauffman and Levin \cite{KL1987} and many others, as mentioned below. Here, a given finite graph represents the possible types of individuals 
(vertices) together with their possibilities of mutation (edges). A fixed, but possibly random, fitness landscape assigns real numbers to the vertices of the graph. 
The evolution of the population is modelled as a random walk on the graph that moves towards higher fitnesses. This can be interpreted as the adaptation of the population to its environment. In contrast to the adaptive dynamics context, this fitness landscape is not dependent on the current state of the population. Adaptive walks move along edges towards neighbours of increasing fitness, according to some transition law, towards a local or global maximum.
In particular, in such models it is not possible for a population to cross a fitness valley. 
This is partially solved by a variation of this model, called \emph{adaptive flight} \cite{NK2011}. It consists in a walk jumping between local fitness maxima, before eventually reaching a global maximum. The questions of the distribution of maxima \cite{NovKru2015}, the typical length of a walk \cite{Orr2003}, or the typical accessibility properties of the fitness landscape \cite{Krug-rev-access-2019,SchKru2014,BerBruShi2016} have been studied under different assumptions on the graph structure, the fitness law, or the transition law of the walk. Moreover, comparisons of these models with actual empirical fitness landscapes have been performed in \cite{Krug-rev-emp-2013}.
As Kraut and Bovier showed \cite{BovKra2018}, adaptive walks and flights arise as the limit of individual-based models of adaptive dynamics, when the large population followed by the rare mutations limit is taken. They also conjecture, and this will be proved in the present article, that similar results hold in the stochastic setting under the mutation rate \eqref{mutation-rate}, as we detail below. \\

In this paper, we consider an individual-based Markov process that models the evolution of a haploid, asexually reproducing population. 
The space of possible traits is given by the vertices of a (possibly directed) finite 
graph $\GG=(V,E)$. 
The evolution of the population is driven by births, deaths, and competition rates, which are fixed and depend on the traits, as well as mutations towards nearest neighbors in the graph $\GG$.
We start with a macroscopic initial condition (that is to say of order $K$, see Definition \ref{def_not}) and we are interested in the stochastic 
process given by the large population limit under the mutation rate \eqref{mutation-rate}.
We describe the time evolution of the orders of magnitude of each sub-population on the $\log K$ time scale, as $K$ tends to infinity. 
We show that the limiting process is deterministic, given by piecewise affine continuous functions, which are determined by an algorithm describing the changes 
in the fitness landscape due to the succession of new resident or emergent types.

This work constitutes an extension of the paper by Kraut and Bovier \cite{BovKra2018} to the stochastic setting. 
They consider the deterministic system resulting from the large population limit of the individual-based model ($K\to\infty$), and let the mutation probability $\mu$ tend to zero. By rescaling the time by $\log(1/\mu)$, they prove that the limiting process is a deterministic adaptive walk that jumps between different equilibria of coexisting traits. A corollary of our results gives the same behaviour, on the $\log K$ time scale, for the stochastic process under the scaling \eqref{mutation-rate} for $\alpha$ larger than the diameter of the graph $\GG$.
Kraut and Bovier also study a variation of the model, where they modify the deterministic system such that the subpopulations can only reproduce when their size lies above a certain threshold $\mu^\alpha$. This limits the radius in which a resident population can foster mutants, and mimics the scaling \eqref{mutation-rate} that we consider. The resulting limiting processes are adaptive flights (which are not restricted to jumping to nearest neighbours), and thus can cross valleys in the fitness landscape and reach a global fitness maximum. We obtain the same behaviour, on the $\log K$ time scale, for the stochastic process under the scaling \eqref{mutation-rate} without any restriction on $\alpha$.

The results of the present paper can also be seen as a generalisation of Theorem 3.2 in \cite{BovCoqSma2018} by Bovier, Coquille and Smadi to any finite trait space. Indeed, 
they consider the graph with vertices $V=\{0,\ldots,L\}$ embedded in $\N$ and choose parameters such that the induced fitness landscape exhibits a valley: mutant
individuals with negative fitness have to be created in order for the population
to reach a trait with positive fitness. Several speeds of the mutation rate are considered, and in particular, when $\alpha>L$, the exit time of the valley is computed on the $\log K$ time scale. This becomes a corollary of our results, and we can give an algorithmic description of the rescaled process for more general graphs endowed with a fitness valley, as we discuss in several examples in Section \ref{sec-examples}.

Our proof heavily relies on couplings of the original process with logistic birth and death processes with non-constant immigration, and the analysis of the later simpler processes on the $\log K$ time scale. This approach was developed by Champagnat, M\'el\'eard and Tran in \cite{htrans19}. They consider an individual-based model for the evolution of a discrete population performing horizontal gene transfer 
and mutations on $V=[0,4]\cap\delta\N,\delta>0$. Their goal is to analyze the trade-off between natural selection, which drives the population to higher birth-rates, and transfer, which drives the population to lower ones. 
Under the mutation rate \eqref{mutation-rate}, they exhibit parameter regimes where different evolutionary outcomes appear, in particular evolutionary suicide and emergence of a cyclic behavior.  As in the present paper, their results characterize the time 
evolution of the orders of magnitude of each sub-population on the $\log K$ time scale, which are shown to be piecewise affine continuous functions whose slopes are given by an algorithm describing the succession of phases when a given type is dominant or resident.
Their proofs provide us with the main ingredients needed for our results. However, the graph structure they choose simplifies the inductions and we 
have to generalise their approach to treat the case of more general graphs, in the proof spirit of Kraut and Bovier \cite{BovKra2018}.

Our results are general, and could be applied to have a better understanding of evolutionary trajectories in complex fitness landscapes. There are now
more and more empirical studies of fitness landscapes (see \cite{de2014empirical} for a comprehensive review of data and tools up to 2014 for instance), and the probability 
and effect of specific mutations in given landscapes are better and better understood.
For instance oriented mutation graphs can stem from mutation bias, through codon usage bias or similar molecular phenomena which make some mutations more probable than others \cite{PloKud}.

We present a series of specific examples where surprising phenomena arise from the geometry of the graph $\GG$ and/or the rate of mutations \eqref{mutation-rate}. Most of them could not happen under a different scaling of mutation rates.

 \begin{itemize}
\item  In Example \ref{example-loop}, we describe a scenario where the ancestry of the resident population consists, with high probability, of back mutations 
towards a previously extinct trait, 
although the mutations that happen in between are not deleterious. In other words, the final resident individuals, say of trait $v$, although they can 
be produced from a wild type directly, come with high probability from a sequence of \emph{non deleterious} mutations which went back to the wild type before mutating 
to $v$. This phenomenon can also happen in the regime \eqref{mut-champagnat}, that is for $\alpha\in(0,1)$, on the mutation time scale ($1/K\mu_K\gg\log K$), where invading mutants fully replace the resident population before a new mutant arises. We show that it can still occur for higher mutation rates of the form \eqref{mutation-rate}, on a $\log K$ time scale, when parameters are chosen such that temporary extinction of the original trait is likely.
Such mutational reversions have been observed (see \cite{DePristo2007} for instance).

\item If evolution and mutation time scales are separated (i.e. in the regime \eqref{mut-champagnat}), mutations occur one at a time, and the number of successive resident traits from the wild type to the type gathering $k$ successively beneficial mutations is $k$. This is not the case if mutations are faster, in which case it is possible to observe either more or less successive resident traits. We will show this in Examples \ref{example-z} and \ref{example2-z}. 

\item In Example \ref{example-anarchy}, we show that adding a new possible mutation path towards a fit trait can increase the time until it appears macroscopically. This is in the spirit of the paradox called \emph{price of anarchy} in game theory or more specifically \emph{Braess paradox} in the study of traffic networks congestion.
 Motter showed that this paradox may often occur in biological and ecological systems \cite{Motter2010}. He studies the removal of part of a metabolic network 
 to ensure its long term persistance, with applications to cancer, antibiotics and metabolic diseases. Another field of application is the food webs management, where selective removal of some species from the network can potentially have a positive outcome of preventing a series of further extinctions \cite{SahMot2011}.

\item Another counter-intuitive phenomenon arising from the mutation rate \eqref{mutation-rate}, presented in Example \ref{example-cycle}, is the possibility 
to observe, for a cyclic clockwise oriented mutation graph, successive counter-clockwise resident populations. 
This means that the macroscopic succession of resident traits is not necessarily representative of the mutation graph. In particular, 
this may call into question the interpretation in terms of 
mutation graphs of some experiments in experimental evolution (see \cite{maddamsetti2015adaptation} for instance).

\item 
In Examples \ref{example-large-jumps} and \ref{example-large-jumps-5sites}, we show that the mutation rate \eqref{mutation-rate} does not restrict the range 
of the corresponding adaptive flights on the trait space, i.e.\ the distance that the limiting process can jump, to $\lfloor \alpha\rfloor$.

\item We finally study the framework of fitness valley crossings. Combining our results with Theorem 3.3 of \cite{BovCoqSma2018}, we construct Examples 
\ref{example-rw2} and \ref{example-rw3}, where effective random walks on the trait space appear on the time scale $K^\beta$, for some positive $\beta$. 
Those limiting adaptive flights arise as a result of a "fast" equilibration on the $\log K$ time scale followed by exponential waiting times until fitness 
valleys get crossed. This makes sense biologically, since there may be traits with positive invasion fitness
that can be reached through several consecutive mutation steps \cite{Len2003,Cow2006}.
 \end{itemize}

The remainder of this paper is organised as follows. 
In Section \ref{sec-results} we define the model and present our results.
In Section \ref{sec-examples} we illustrate our results by a series of examples describing surprising phenomena arising from the geometry and/or the rate of mutations.
Section \ref{sec-proof} is devoted to the proofs.
In the Appendix, we present and extend some technical results.


\section{Convergence on the $\log K$-time scale}\label{sec-results}
\subsection{Model}
We consider an individual-based Markov process that models the evolution of a haploid, asexually reproducing population. 
The space of possible traits is given by the vertices of a (possibly directed) finite 
graph $\GG=(V,E)$. 

For all traits $v,w\in V$ and every $K\in\N$, we introduce the following parameters:
\begin{itemize}
\item $b_v\in\R_+$, the birth rate of an individual of trait $v$,
\item $d_v\in\R_+$, the (natural) death rate of an individual of trait $v$,
\item $c^K_{v,w}=c_{v,w}/K\in\R_+$, the competition imposed by an individual of trait $w$  onto an individual of trait $v$,
\item $\mu_K\in[0,1]$, the probability of mutation at a birth event,
\item $m(v,\cdot)\in\MM_p(V)$, the law of the trait of a mutant offspring produced by an individual of trait $v$.
\end{itemize}

The process $N^K$ describes the state of the population, where $N^K_v(t)$ denotes the number of individuals of 
trait $v\in V$ alive at time $t\geq 0$.
We assume that edges in $E$ mark the possibility of mutation and hence $m(v,w)>0$ if and only if $(v,w)\in E$.

\begin{remark}
We could also allow for $\mu_K$ to depend on $v\in V$ as long as $\mu_K(v)=\mu_Kh(v)$ for some strictly positive function $h$ that is independent of $K$. However, this would not change the characterisation of the limit, and hence we assume a constant $\mu_K$ to simplify the notation.
\end{remark}

Moreover, we assume that, for every $v\in V$, $c_{v,v}>0$. The parameter $K$ is scaling the competitive pressure and, 
through this self-competition, fixes the equilibrium size of the population to the order of $K$. $K$ is sometimes called \textit{carrying capacity} and can be interpreted as a scaling parameter for the available sources of food or space.

As a consequence of our parameter definitions, the process $N^K$ is characterised by its infinitesimal generator:
\begin{align}\label{generator}
\LL^K\phi(N)=&\sum_{v\in V}(\phi(N+\delta_v)-\phi(N))\left(N_vb_v(1-\mu_K)+\sum_{w\in V}N_wb_w\mu_Km(w,v)\right)\notag\\
&+\sum_{v\in V}(\phi(N-\delta_v)-\phi(N))N_v\left(d_v+\sum_{w\in V}c^K_{v,w}N_w\right),
\end{align}
where $\phi:\N^V\to\R$ is measurable and bounded. Such processes have been explicitely constructed in terms of Poisson random measures in \cite{FM04}.

Due to the scaling of the competition $c^K$, the equilibrium population is of order $K$.
Since the mutation probability $\mu_K$ tends to zero as $K\to\infty$, the process $N^K/K$ converges (on finite time intervals)  to the mutation-free Lotka-Volterra system \eqref{LV_system} involving all initial coexisting resident traits. 
We are interested in the long-term evolution of the population and want to study successive invasions by new mutant populations. Given the fact that a mutant population that is initially of order $K^\gamma$, $\gamma<1$, needs a time of order $\log K$ to grow exponentially to the order of $K$, we have to rescale the time by $\log K$ to obtain a non trivial limit.

It is convenient to describe the population size of a certain trait $v\in V$ by its $K$-exponent
\begin{align}
\beta^K_v(t):=\frac{\log(1+N^K_v(t\log K))}{\log K},
\end{align}
which is equivalent to $N^K_v(t\log K)=K^{\beta^K_v(t)}-1$. Since the population size is restricted to order $K$ by the competition, $\beta^K_v$ 
ranges between $0$ and $1$, as $K\to\infty$ (see Corollary \ref{betabound} for a rigorous statement).

For the sake of readability, we now introduce the terminology we will use in the sequel.

\begin{definition} \label{def_not} \emph{}
	\begin{enumerate}
		\item A trait $v\in V$ with exponent $\beta^K_v$ is called \emph{macroscopic} if, for every $\eps>0$, there exists $K_\eps$ such that, 
		for every $K\geq K_\eps$, $\beta^K_v>1-\eps$ .
		\item A trait that is not macroscopic is called \emph{microscopic}.
		\item The set of \emph{living traits} is the set $\{v\in V : \beta^K_v>0\}$.
	\end{enumerate}
\end{definition}

When $K$ is large enough, the macroscopic traits interact on any finite time interval according to the corresponding mutation-free Lotka-Volterra 
system (see Chapter 11, Theorem 2.1 in \cite{ethiermarkov} for the proof of this law of large numbers): Let $\v \subset V$, 
then the mutation-free Lotka-Volterra system associated to $\v$ is 
\begin{equation} \label{LV_system}
 \dot{n}_w(t)= \left( b_w- d_w - \sum_{v \in \v}c_{w,v}n_v(t) \right) n_w(t), \quad w \in \v, \ t \geq 0. 
\end{equation}

For a subset $\v\subset V$ of traits,
we denote by $\bar{n}(\v)\in\R_+^V$ the unique equilibrium of the Lotka-Volterra system \eqref{LV_system}, when it exists, and where to simplify notations, we extend it by
$\bar{n}_w(\v)=0$ for $w \notin \v$.
In the case where $\v=\{v\}$, we obtain from classical results on Lotka-Volterra models (see \cite{Cha06} for instance) 
$$\bar{n}_v(v)=(b_v-d_v)/c_{v,v} \vee 0.$$
If $\v$ denotes the set of macroscopic traits, we call the traits $v\in\v$ such that $\bar{n}_v(\v)>0$ \textit{resident}.

The approximate rate at which a mutant of trait $w$ grows in a population of coexisting resident traits $\v$ is called \textit{invasion fitness} and is denoted by $f_{w,\v}$, where
\begin{align*}
f_{w,\v}:=b_w-d_w-\sum_{v\in\v}c_{w,v}\bar{n}_v(\v).
\end{align*}
If $f_{w,\v}>0$, the trait $w$ is called \emph{fit}.
If $f_{w,\v}<0$, the trait $w$ is called \emph{unfit}. 
The case $f_{w,\v}=0$ will be excluded (see Remark \ref{rmk_termination}).

Mutants can be produced along (directed) edges of the graph. We denote by $d(v,w)$ the graph distance, i.e.\  the length of the shortest (directed) path 
from $v$ to $w$ in $\GG=(V,E)$. For a subset $\v\subset V$ we define 
$$d(\v,w):=\min_{v\in\v}d(v,w) \quad \text{and} \quad d(w,\v):=\min_{v\in \v}d(w,v).$$

\subsection{Results}

Let a finite graph $\GG=(V,E)$ be given and assume that $\alpha \in \R_+^* \setminus \N$ and $f_{w,\v} \neq 0$ for any $w \in V$ and $\v \subset V$
(see Remark \ref{rmk_termination}). 
The two following results concern the convergence of the orders of the different subpopulation sizes to a piecewise linear trajectory, whose slopes and times of slope changes can be explicitely expressed in terms of the parameters. 

\begin{theorem}\label{ThmConv}
Let a finite graph $\GG=(V,E)$ and $\alpha \in \R_+^* \setminus \N$ be given and consider the model defined by \eqref{generator}. Assume that
$f_{w,\v} \neq 0$ for any $w \in V$ and $\v \subset V$. 
Let $\v_0\subset V$ and assume that, for every $w\in V$,
\begin{align}\label{cond-init}
\beta^K_w(0)\rightarrow\left(1-\frac{d(\v_0,w)}{\alpha}\right)_+, \quad (K\to\infty) \quad \text{in probability}.
\end{align}
Then, for all $T>0$, as $K\to\infty$, the sequence $((\beta^K_w(t),w\in V),t\in[0,T\land T_0])$ converges in probability in 
$\mathbb{D}([0,T\land T_0],\R_+^V)$ to a deterministic, piecewise affine, continuous function $((\beta_w(t),w\in V),t\in[0,T\land T_0])$, which is defined as follows:
\begin{itemize}
\item[(i)] If the mutation-free Lotka-Volterra system \eqref{LV_system} associated to $\v_0$ has a unique positive globally attractive equilibrium, 
the initial condition of $\beta$ is set to $\beta_w(0):=\left(1-\frac{d(\v_0,w)}{\alpha}\right)_+$.
Otherwise, the construction is stopped and $T_0$ is set to $0$.
\item[(ii)] The increasing sequence of invasion times is denoted by $(s_k)_{k\geq0}$, where $s_0:=0$ and, for $k\geq1$,
\begin{align*}
s_k:=\inf\{t>s_{k-1}:\exists w\in V\backslash\v_{k-1}:\beta_w(t)=1\}.
\end{align*}
Here, $\v_k$ denotes the set of coexisting resident traits of the Lotka-Volterra system that includes $\v_{k-1}$ and the trait $w\in V\backslash\v_{k-1}$ that satisfies $\beta_w(s_k)=1$.
\item[(iii)] For $s_{k-1}\leq t\leq s_k$, for any $w\in V$, $\beta_w(t)$ is defined by
\begin{align}\label{beta}
\beta_w(t):=\max_{\substack{u\in V}}\left[\beta_u(s_{k-1})+(t-t_{u,k}\land t)f_{u,\v_{k-1}}-\frac{d(u,w)}{\alpha}\right]\lor 0,
\end{align}
where, for any $w \in V$,
\begin{align} \label{deftwk}
t_{w,k}:=\begin{cases}\inf\{t\geq s_{k-1}:\exists\ u\in V: d(u,w)=1, \beta_u(t)=\frac{1}{\alpha}\}
&\text{if }\beta_w(s_{k-1})=0\\s_{k-1}&\text{else}\end{cases}
\end{align}
is the first time in $[s_{k-1},s_k]$ when this trait arise.
\item[(iv)] The inductive construction is stopped and $T_0$ is set to $s_k$ if
\begin{itemize}
\item[(a)] there is more than one $w\in V\backslash\v_{k-1}$ such that $\beta_w(s_k)=1$;
\item[(b)] the Lotka-Volterra system including $\v_{k-1}$ and the unique 
$w\in V\backslash\v_{k-1}$ such that $\beta_w(s_k)=1$ does not have a unique stable equilibrium;
\item[(c)] there exists $w\in V\backslash\v_{k-1}$ such that $\beta_w(s_k)=0$ and 
$\beta_w(s_k-\e)>0$ for all $\e>0$ small enough. 
\item[(d)] there exists $w\in V\backslash\v_{k-1}$ such that $s_k=t_{w,k}$.
\end{itemize}
\end{itemize}
\end{theorem}

\begin{remark}\label{rmk_termination}
Notice that conditions $(a)$, $(c)$, and $(d)$ of point $(iv)$ are here to exclude very specific and non generic cases where one coordinate reaches $1$ 
while another reaches $1$ or reaches $0$ from above, or a new trait arises at the exact same time. They are difficult to handle for technical reasons.

Moreover, we exclude the cases where $\alpha\in\N$. They would produce mutant populations, at distance $\alpha$ from the resident traits, 
that can neither be approximated by sub- nor super-critical branching processes. The same applies to the case $f_{w,\v}=0$, where the population can both grow 
and shrink due to fluctuations.
\end{remark}

\begin{remark}
The $t_{w,k}$ do not keep track of traits that die out in $[s_{k-1},s_k]$ and then reappear. However, since the fitnesses do 
not change between invasions, such a trait would have {a} negative invasion fitness (else it would not die out). 
Hence, it would not start growing on its own if it reappears, but only follow along another trait due to mutants. 
It would therefore not contribute to the maximum over $u\in V$ in \eqref{beta}.
\end{remark}

\begin{proposition} \label{CorEquilibria}
Under the same assumptions and with the same notations as in Theorem \ref{ThmConv}, 
for all $T>0$, as $K\to\infty$, the sequence $((N^K_w(t\log K)/K,w\in V),t\in[0,T\land T_0])$ 
converges in probability in $\mathbb{D}([0,T\land T_0]\backslash\{s_k,k\geq1\}, \R_+^V)$ 
to a deterministic jump process $((N_w(t),w\in V),t\in[0,T\land T_0])$, which is defined as follows:
\begin{itemize}
\item[(i)] For $t\in[0,T_0]$, $N(t)$ jumps between different Lotka-Volterra equilibria according to
\begin{align*}
N_w(t):=\sum_{k\in\N:s_{k+1}\leq T_0}\mathbf{1}_{s_k\leq t<s_{k+1}}\mathbf{1}_{w\in\v_k}\bar{n}_w(\v_k).
\end{align*}
\item[(ii)]  The invasion times $s_k$ and the times $t_{w,k}$ when new mutants arise 
can be calculated as follows.
We define the increasing sequence $(\tau_\ell,{\ell\geq0})=\{s_k,{k\geq0}\}\cup\{t_{w,k},{w\in V,k\geq0}\}$ of invasion times or appearance times of new mutants, 
and $(M_\ell, \ell\geq0)$ the sets of living traits in the time interval $(\tau_\ell,\tau_{\ell+1}]$. Initially, $\tau_0=s_0=0$ 
and, according to \eqref{cond-init}, $M_0=\{w\in V:d(\v_0,w)<\alpha\}=\{w\in V:\beta_w(0)>0\}$. For $s_{k-1}\leq\tau_{\ell-1}<s_k$, $\tau_\ell$ is defined as
\begin{align*}
\tau_\ell:=s_k\land\min\{t_{w,k}:w\in V, t_{w,k}>\tau_{l-1}\}.
\end{align*}

Given $\tau_\ell$ and $M_{\ell-1}$, we set $M_\ell:=(M_{\ell-1}\backslash\{w\in V:\beta_w(\tau_\ell)=0\})\cup\{w\in V:\tau_\ell=t_{w,k}\}$. $\tau_\ell$ is then given by
\begin{align}\label{Dtau-formula}
\tau_\ell-\tau_{\ell-1}=\min_{\substack{w\in M_{\ell-1}:\\f_{w,\v_{\ell-1}}>0}}
\frac{\left(1\land\frac{d(w,V\backslash M_{\ell-1})}{\alpha}\right)-\beta_w(\tau_{\ell-1})}{f_{w,\v_{\ell-1}}}.
\end{align}
\end{itemize}
\end{proposition}

\begin{remark} We could allow for more general initial conditions of the form
\begin{align*}
\beta^K_w(0)\to\tilde{\beta}_w\in[0,1],
\end{align*}
with $\tilde{\beta}_w$, $w\in V$, deterministic and $\v_0:=\{w\in V:\tilde{\beta}_w=1\}\neq\emptyset$. 
An inductive application of Lemma \ref{lemmaB4}, similar to the induction proving \eqref{inductiontoprove}, implies that within a 
time of order 1, for all $w\in V$, $\beta^K_w\cong \max_{u\in V}\{\tilde{\beta_u}-d(u,w)/\alpha$\}. 
We therefore set $\beta_w(0):=\max_{u\in V}\{\tilde{\beta_u}-d(u,w)/\alpha\}$ in Theorem \ref{ThmConv} and $M_0:=\{w\in V:\exists\ u\in V\text{ s.t.\ }\beta_u(0)>0, 
\beta_u(0)\geq d(u,w)/\alpha\}$ in Proposition \ref{CorEquilibria}. The rest of the results remains unchanged.
\end{remark}

\begin{remark}
The limiting jump process $N(t)$ resembles an adaptive walk or flight, as studied in \cite{Orr2003,NK2011,SchKru2014,NovKru2015,BerBruShi2016}. For a constant competition kernel $c_{v,w}\equiv c$, we consider the fixed fitness landscape given by $r_v=b_v-d_v$. Since in this case $f_{w,v}=r_w-r_v$, the process jumps along edges towards traits of increasing fitness $r$.
\end{remark}

The above results are in the vein of Theorem 2.1 and Corollary 2.3 in \cite{htrans19}. There are however many differences between the setting considered in \cite{htrans19} and our setting.

Due to the horizontal transfer between individuals, Champagnat and coauthors obtained trajectories where a "dominant" population, i.e.\ with the size 
of highest order, could be non resident, i.e.\ of order negligible with respect to $K$. They could also witness extinction on a 
$\log K$ time scale as well as evolutionary 
suicide. The absence of horizontal transfer in our case prevents such behaviours.

We consider a general finite graph of mutations with 
possible back mutations, whereas their graph was embedded in $\Z$ and did not allow for back mutations. We also allow for the coexistence
of several resident traits in the population at equilibrium. The two main difficulties in the proofs compared to \cite{htrans19} are thus to handle the generality 
of the graph of mutations, and to extend some approximation results to the multidimensional case.


\section{Surprising phenomena arising from geometry and mutation rate}
 \label{sec-examples}
 
In this section, we present some non intuitive behaviours of the population process, which stem from the mutation scale or the generality of the 
mutational graph that we allow for. They are direct applications of Theorem \ref{ThmConv} and Proposition \ref{CorEquilibria}, and provide explicit computations of exponents \eqref{beta} and time intervals \eqref{Dtau-formula}.

Several examples are build on directed graphs. Although this is not a necessary condition to obtain the desired phenomena, it allows a simplified study (especially of the decay phases).

We first introduce some notations for the sake of readability.

\begin{definition} Let $w,v\in V$ and $\v\subset V$. We write
	\begin{enumerate}
		\item \emph{with high probability} to mean "with a probability converging to 1 as $K\to\infty$",
		\item $w>\v$ if and only if $f_{w,\v}>0$, that is if $w$ can invade in $\v$,
		\item $w<\v$ if and only if $f_{w,\v}<0$, that is if $w$ cannot invade in $\v$,
		\item $w\gg v$ (or $v\ll w$) if and only if $f_{w,v}>0$ and $f_{v,w}<0$, that is if $w$ can invade in $v$ and fixate,
		\item $w\equiv v$ if and only if $f_{w,v}>0$ and $f_{v,w}>0$, that is if $w$ and $v$ can coexist,
		\item $w \frown v$ if and only if $f_{w,v}<0$ and $f_{v,w}<0$, that is if $w$ and $v$ can neither invade in each other.
	\end{enumerate}
\end{definition}

\subsection{Back mutations before adaptation}
In the following, we build an example where the ancestry of the resident population comes from back mutations from an ancestral trait, even if the mutations happening 
in between are not deleterious.

\begin{figure}
	\centering
	\hspace{1cm} \includegraphics[width=.2\textwidth]{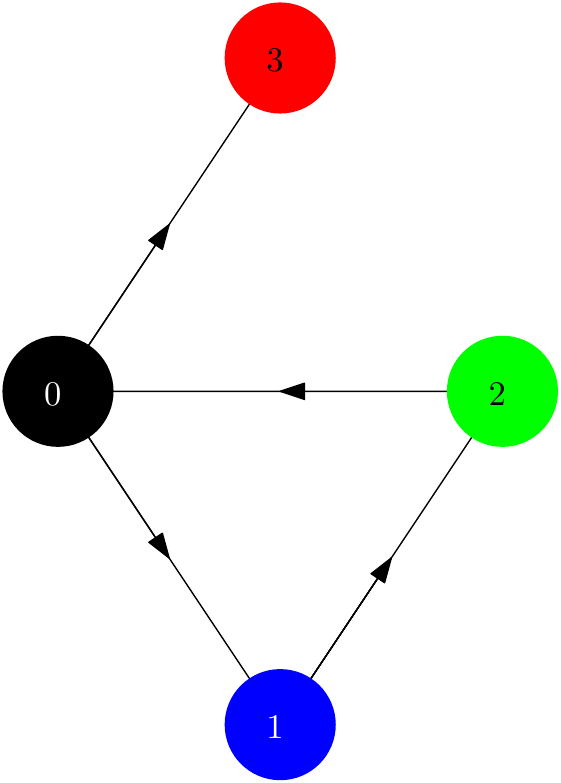} \hspace{1cm}
	\includegraphics[width=.25\textwidth]{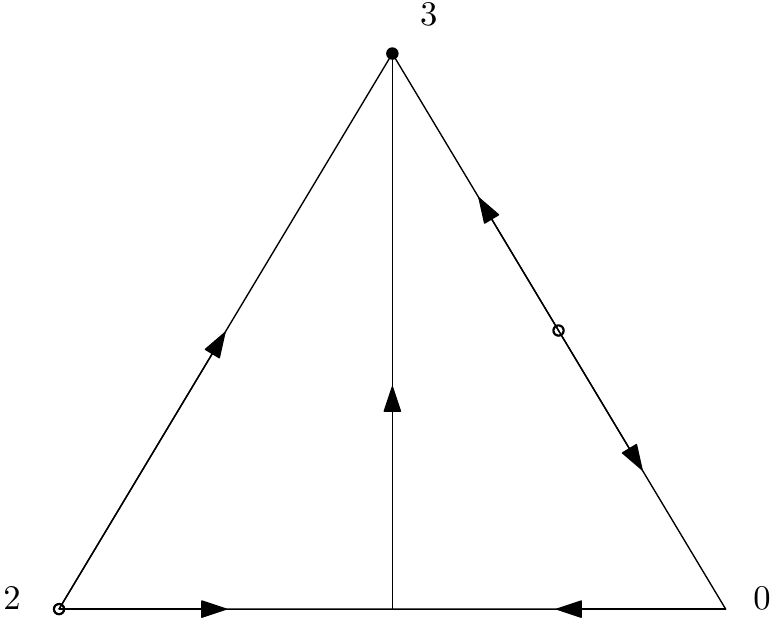}
	\includegraphics[width=.8\textwidth]{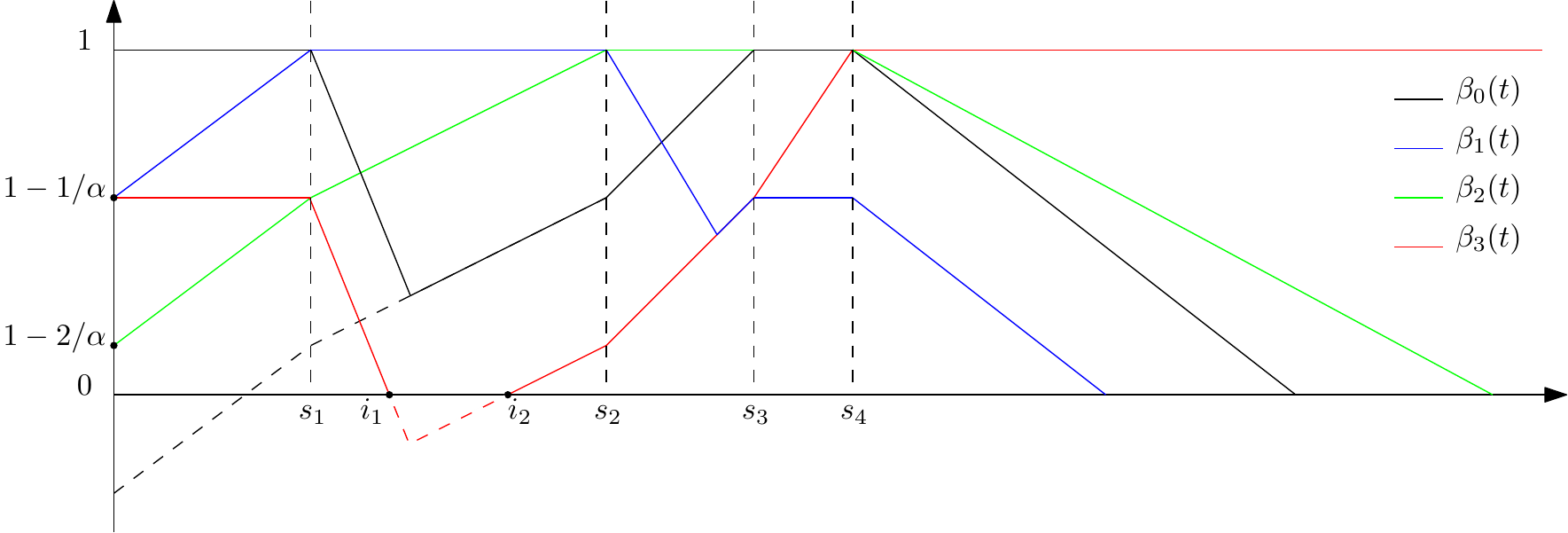}
	\caption
	{Graph $\GG$, phase portrait of traits 0-2-3, and exponents $\beta(t)$ of Example \ref{example-loop}.}
	\label{fig-loop}
\end{figure}

\begin{example}\label{example-loop}
	Let us consider the graph $\GG$ depicted on Figure \ref{fig-loop} where $V=\{0,1,2,3\}$ and $E=\{[0,1],[1,2],[2,0],[0,3]\}$.
	Let $\alpha>2$, an initial condition given by $(\bar{n}(0),0,0,0)$ and a fitness landscape given by
	\begin{align}
	&{0\ll1\ll2, \quad 3\frown 0,\quad 3\frown 1} \nonumber \\
	&{
	0\equiv 2, \quad
	3>\{0,2\},\quad
	2<3} \label{loop-cond-302}\\
	&f_{2,0}<2f_{1,0}\label{loop-cond-1resident}\\
	&f_{0,1}\geq f_{3,1},\quad f_{2,0}\leq f_{1,0}\label{SNNC1}\\
	&i_1:=\frac{1-1/\alpha}{f_{0,1}}<\frac{-(1-4/\alpha)}{f_{2,1}}=:i_2 \label{loop-cond-whp}
	\end{align}
\end{example}

In this case, Proposition \ref{CorEquilibria} implies that on the $\log K$ time scale, the rescaled macroscopic population 
then jumps from traits $0-1-2$ then to coexistence between $0$ and $2$, followed by the invasion and fixation of $3$ which 
is produced with high probability, due to Condition \eqref{loop-cond-whp}, by individuals of type $0$ which have the sequence 
$0-1-2$ as ancestry. 
In other words, the final resident individuals of trait 3, although they can be produced by 
individuals of trait 0 directly, come from a sequence of mutations which went around the loop $0-1-2$ of $\GG$.
Conditions \eqref{loop-cond-302}, summarized on Figure \ref{fig-loop},  imply phase portrait number 8 in the classification of Zeeman \cite{Zee1993}. 
Condition \eqref{loop-cond-1resident} ensures that trait 1 becomes resident before 2.
Condition \eqref{SNNC1} is not necessary but allows to simplify the setting.
The exponents are drawn on Figure \ref{fig-loop}.

Note that this scenario can also happen in the rare mutation regime considered in \cite{Cha06} (for example $\alpha\in(0,1)$): 
the average waiting time until a mutant of type 1 appears is then of order $O(1/K\mu_K)=O(K^{-1+1/\alpha})\gg\log K$. Once 
it has appeared, it survives with positive probability and the succession of invasions and fixations above takes place on 
the $\log K$ time scale, separated by mutation events on the $K^{-1+1/\alpha}$ time scale. 
What is new in our case is that such a scenario 
can still take place for higher mutation rates than the ones considered in \cite{Cha06}, and on a $\log K$ time scale.

\subsection{Non-intuitive mutational pathways in the high mutation framework}
\subsubsection{Longer or shorter path than expected}
If evolution and mutation time scales are separated (i.e. in the rare mutation regime), mutations occur one at a time, 
and the number of successive resident traits from the wild type to the type gathering $k$ successively beneficial mutations 
is $k$. This is not the case if mutations are faster, in which case it is possible to observe either more or less resident traits, as the following examples show.

\begin{figure}
	\centering
	\includegraphics[width=.2\textwidth]{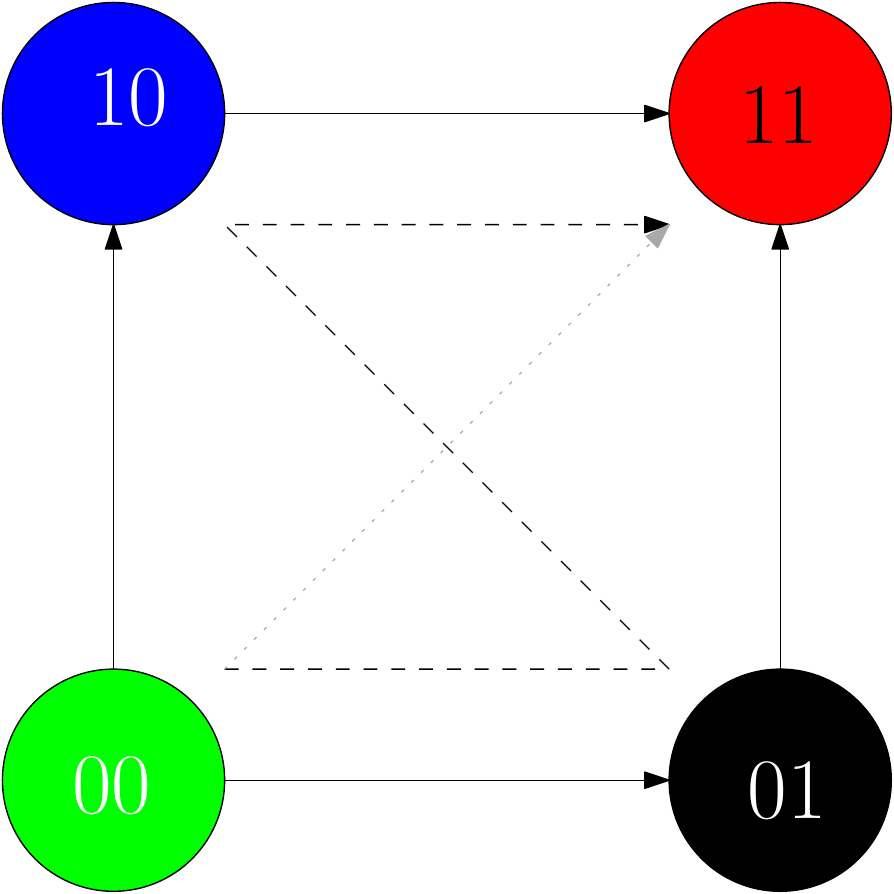}\\
	\includegraphics[width=.7\textwidth]{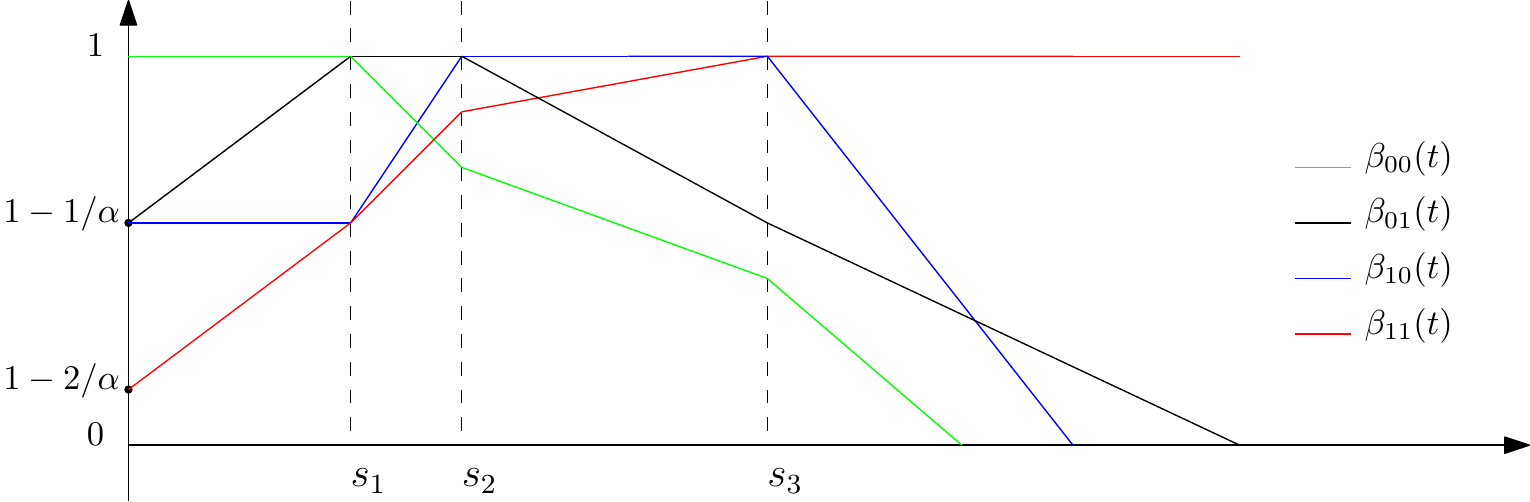}
	\includegraphics[width=.7\textwidth]{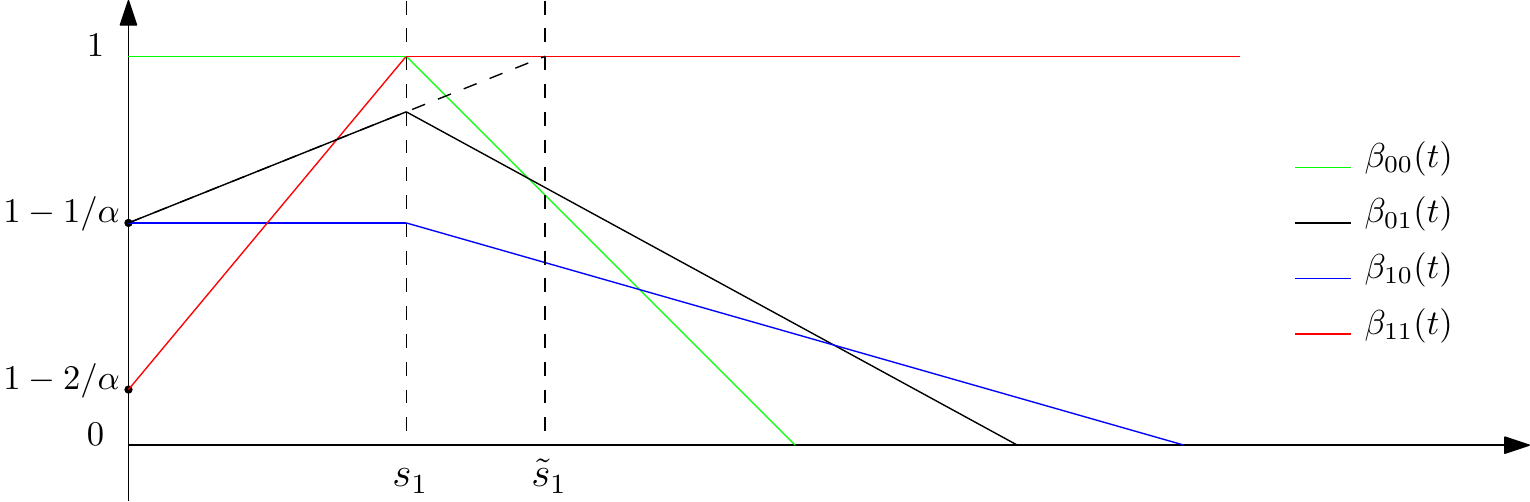}
	\caption
	{Graph $\GG$ and exponents $\beta(t)$ for Examples \ref{example-z} and \ref{example2-z}.}
	\label{fig-z}
\end{figure}

\begin{example} \label{example-z}
	Let us consider the directed graph $\GG$ depicted on Figure \ref{fig-z}, where 
	$V=\{00,01,10,11\}$ and $E=\{[00,01],[00,10],[01,11],[10,11]\}$. 
	Let $\alpha>2$, an initial condition given by $(\bar{n}(00),0,0,0)$ and a fitness landscape given by
	\begin{align}
	&{00\ll01\ll10\ll11},\quad and \quad 01\ll11 \nonumber \\
	&10,11\frown 00\nonumber \\
	&f_{11,00}<f_{01,00}\label{z-cond1}\\
	&f_{10,01}> f_{11,01}\label{z-cond2}
	\end{align}
\end{example}

In this case, in the rare mutation regime, the rescaled macroscopic population jumps along $00-01-11$.

In the regime of Theorem \ref{ThmConv}, Proposition \ref{CorEquilibria} implies that the rescaled macroscopic population jumps along $00-01-10-11$ on the $\log K$ time scale.
More precisely, the exponents are drawn on Figure \ref{fig-z}.
Note that Condition \eqref{z-cond1} ensures that 11 does not invade before 01, it is not necessary but allows to simplify the setting.
Condition  \eqref{z-cond2} ensures that 11 does not invade before 10.

\begin{example} \label{example2-z}
	Let us consider the directed graph $\GG$ depicted on Figure \ref{fig-z}, where $V=\{00,01,10,11\}$ and $E=\{[00,01],[00,10],[01,11],[10,11]\}$. 
	Let $\alpha>2$, an initial condition given by $(\bar{n}(00),0,0,0)$ and a fitness landscape given by
	\begin{align}
	&{01>00,\quad 11\gg00} \nonumber \\
	&{10<00} \label{SNNC2}\\
	&{01,10<11} \nonumber \\
	&\frac{2}{f_{11,00}}<\frac{1}{f_{01,00}}\label{z-cond}
	\end{align}
\end{example}

In this case, in the rare mutation regime, the rescaled macroscopic population still jumps along $00-01-11$, under the additional assumption that $f_{11,01}>0$.

In the regime of Theorem \ref{ThmConv}, Proposition \ref{CorEquilibria} implies that the rescaled macroscopic population directly jumps from $00$ to $11$ on the $\log K$ time scale.
More precisely, the exponents are drawn on Figure \ref{fig-z}. Condition \eqref{z-cond} ensures that $11$ fixates before $01$. 
{Condition \eqref{SNNC2} is not necessary but allows to simplify the setting. Note that equation \eqref{Dtau-formula} implies 
that $s_1={2}/{f_{11,00}}$ and $\tilde s_1={1}/{f_{01,00}}$.

\subsubsection{Price of anarchy}

We build an example where adding a new possible mutation path to a fit trait increases the time until it appears macroscopically.

\begin{example} \label{example-anarchy}
	Let us consider the graph $\GG$ depicted on Figure \ref{fig-anarchy}, where $V=\{1,2a,2b,3\}$ and the edge set is either 
	$E_1=\{[1,2a],[2a,3],[2b,3],[3,2b]\}$ or \\
	$E_2=\{[1,2a],[2a,3],[2b,3],[2a,2b],[3,2b]\}$.
	Let $\alpha>3$, an initial condition given by $(\bar{n}(1),0,0,0)$ and a fitness landscape given by
	\begin{align}
	&1\ll2a\ll3,\quad and \quad 2a\ll2b\\
	&1<2b,\quad and
	,\quad 1, 2b<3 \nonumber \\
	&f_{2a,1}\geq f_{3,1},f_{2b,1},\quad\text{ and }\quad
	\frac{1}{f_{2b,2a}}<\frac{1}{f_{3,2a}}<\frac{2}{f_{2b,2a}}\label{anarchy-cond1}\\
	&{0<}f_{3,2b}<f_{3,2a}.\label{anarchy-cond2}
	\end{align}
\end{example}

In this case, if the edge set is $E_1$, Proposition \ref{CorEquilibria} implies that the rescaled macroscopic population jumps along traits $1-2a-3$ in a time $t_1$ on the $\log K$ time scale.
But if the edge set is $E_2$, the population jumps along $1-2a-2b-3$ and the time to reach $3$ is $t_2>t_1$.
More precisely, the exponents are drawn on Figure \ref{fig-anarchy}. { Condition \eqref{anarchy-cond1} ensures that $2b$ invades first when the edge set is $E_2$ but not when it is $E_1$, in other words $\beta_{2b}$ reaches 1 before $\beta_{3}$ if started at $1-1/\alpha$ but not at $1-2/\alpha$. And Condition \eqref{anarchy-cond2} enlarges the time of fixation of 3. Note that the first inequality in Condition \eqref{anarchy-cond1} is not necessary but allows to simplify the second one.
Moreover, observe that equation \eqref{Dtau-formula} implies $\tilde s_2-\tilde s_1={1}/{f_{2b,2a}}$ and $ s_2- s_1={1}/{f_{3,2a}}$.
Note that in the rare mutation regime we can observe this phenomenon on the mutation time scale, but only with 
probability strictly smaller than 1, since both $2b$ and $3$ are fit with respect to $2a$ and can both invade with positive probability once they are produced.

\begin{figure}
	\centering
	\includegraphics[width=.2\textwidth]{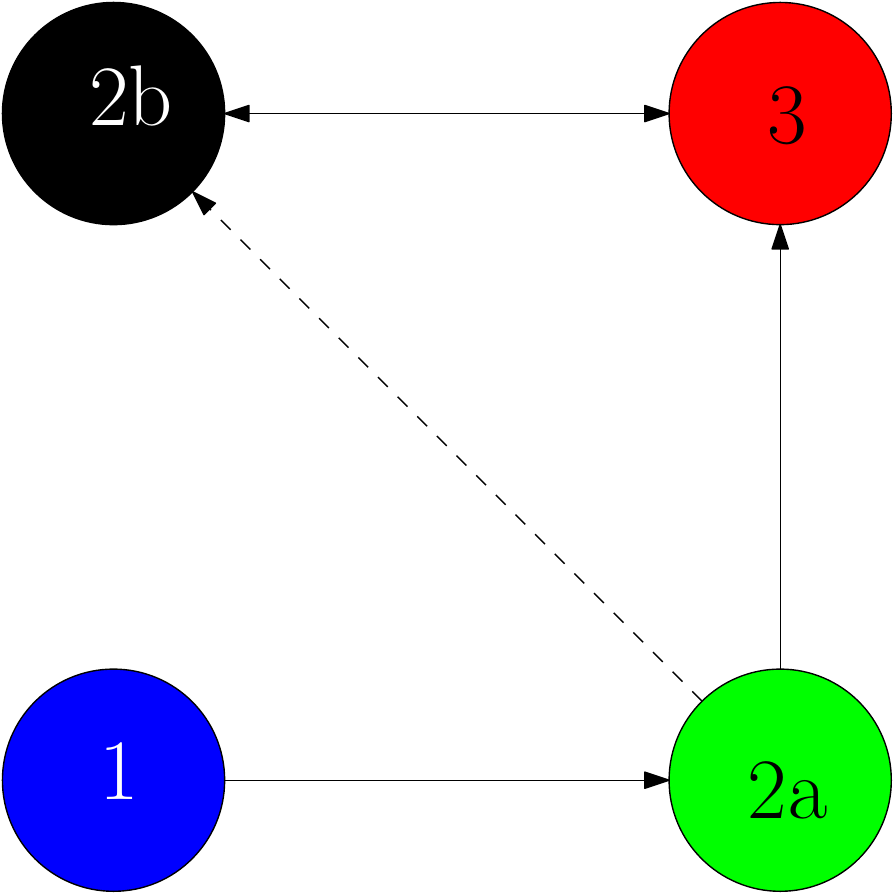}\\
	\includegraphics[width=.7\textwidth]{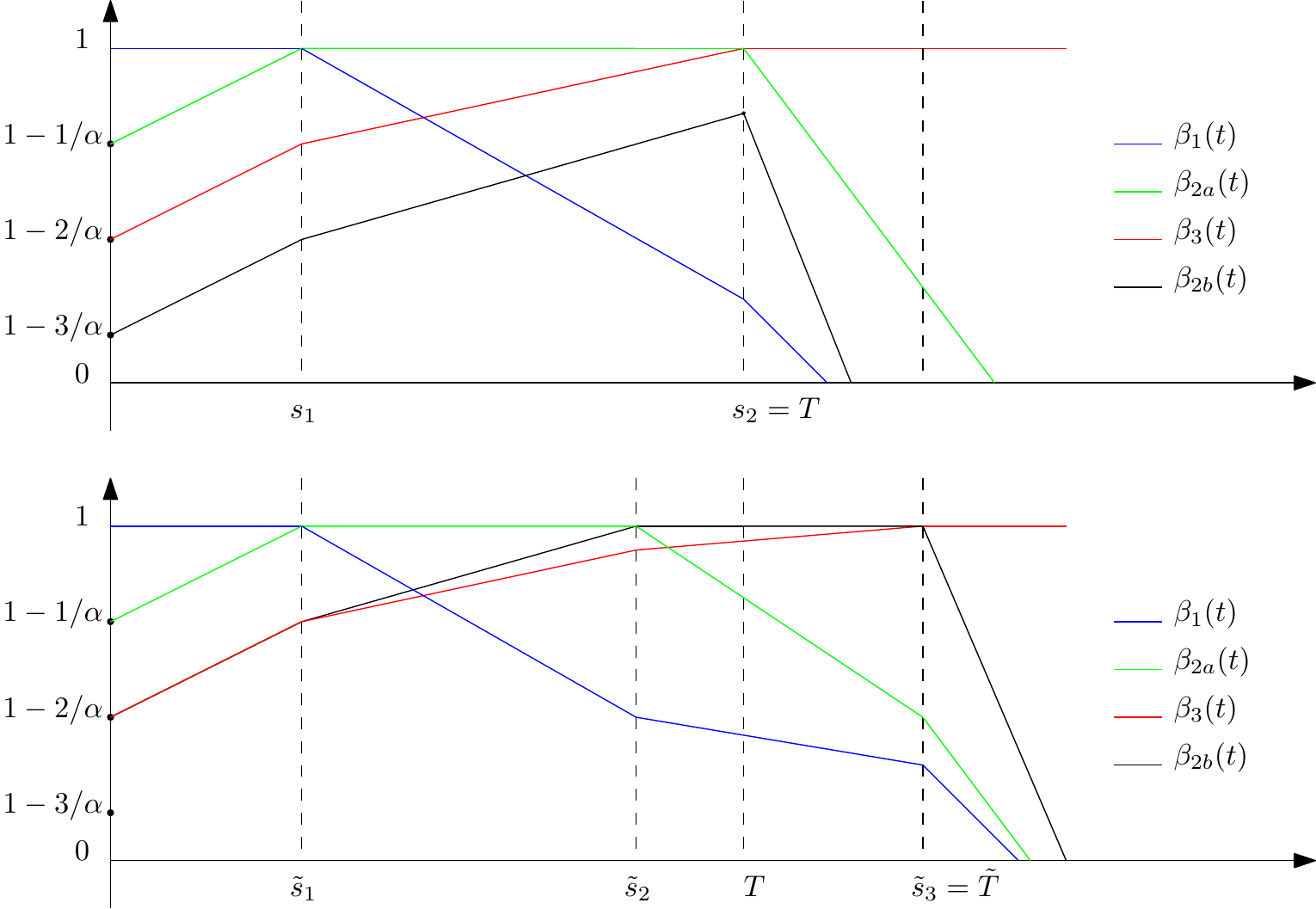}
	\caption
	{Graph $\GG$ and exponents $\beta(t)$ for Example \ref{example-anarchy}, with edge set $E_1$ (above) and $E_2$ (below). }
	\label{fig-anarchy}
\end{figure}

\subsubsection{Counter cycle}

\begin{example} \label{example-cycle}
	Let us consider the graph $\GG$ depicted on Figure \ref{fig-cycle}, where $V=\{1,2,3\}$ and the edge set is $E=\{[1,2],[2,3],[3,1]\}$.
	Let $\alpha>2$, an initial condition given by $(\bar{n}(1),0,0)$ and a fitness landscape given by
	\begin{align*}
	{1\gg 2,\quad 2\gg3 ,\quad 3\gg1.}
	\end{align*}
\end{example}

In this case, Proposition \ref{CorEquilibria} implies that the rescaled macroscopic population jumps along traits $1-3-2$ (in the clockwise sense) although the mutations are directed counterclockwise.
More precisely, the exponents are drawn on Figure \ref{fig-cycle}. 
Moreover, if Conditions \eqref{cycle-cond} below are fulfilled the period is shorter and shorter, an acceleration takes place, as it is depicted on Figure \ref{fig-cycle}.
	\begin{align}
	f_{2,3}>-f_{1,3}\nonumber\\
	f_{1,2}>-f_{3,2}\nonumber\\
	f_{3,1}>-f_{2,1}\label{cycle-cond}
	\end{align}

Note that in the rare mutation regime, with the chosen parameters, there would be no evolution since $2<1$. Moreover, there are no parameters such that counter cyclic or accelerating behaviors could arise.

\begin{figure}
	\centering
	\includegraphics[width=.2\textwidth]{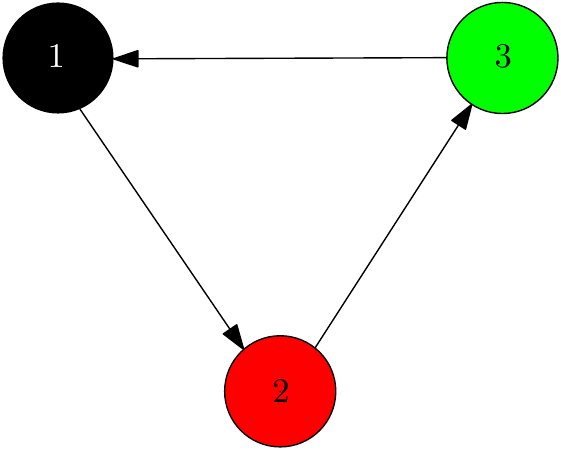}\\
	\includegraphics[width=.8\textwidth]{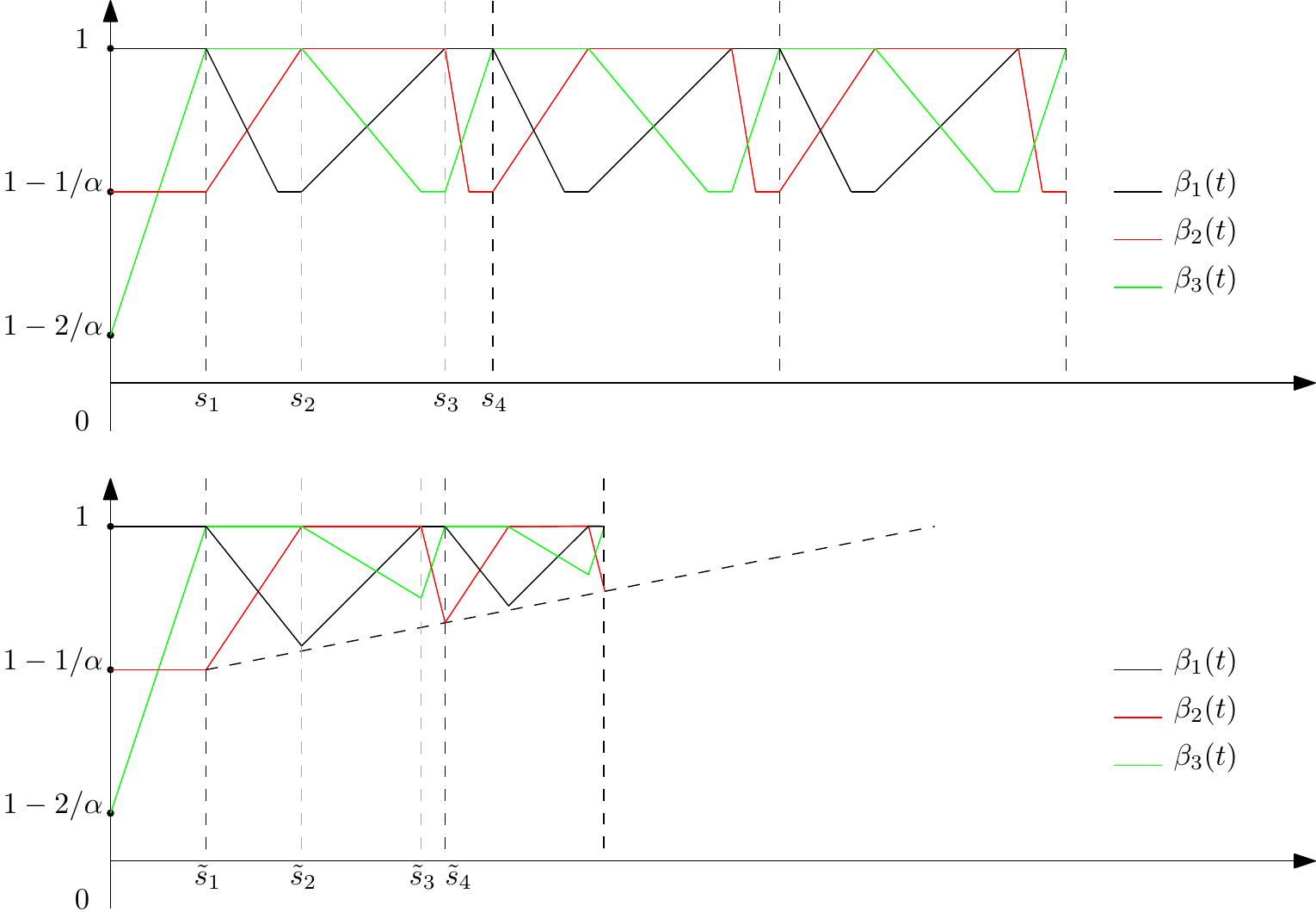}
	\caption
	{Graph $\GG$ and exponents $\beta(t)$ for Example \ref{example-cycle}, without Assumption \ref{cycle-cond} (above) and with Assumption \ref{cycle-cond} (below).}
	\label{fig-cycle}
\end{figure}

\subsection{Arbitrary large jumps on the $\log K$-time scale}
A natural question to ask is if the "cut-off" $\alpha$ restricts the range of the jumps,  on the $\log K$ time-scale, to traits which are at a distance less than $\alpha$. The answer is no, as the following example shows. 

\begin{example} \label{example-large-jumps}
	Let us consider the graph $\GG$ depicted on Figure \ref{fig-large-jumps}, where $V=\{0,1,2,3,4\}$ and $E=\{[0,1],[1,2],[2,3],[3,4]\}$. 
	Let $\alpha\in(3,4)$, an initial condition given by $(\bar{n}(0),0,\ldots,0)$ and a fitness landscape given by
	\begin{align}
	&1,2<0 ,\quad 3,4>0,\quad 0,1,2,3<4 \nonumber \\
	&\frac1{f_{4,0}}+\frac{-1+4/\alpha}{f_{3,0}}<\frac{3/\alpha}{ f_{3,0}} \label{large-jumps-cond}
	\end{align}
\end{example}

In this case, the cut-off is in between traits 3 and 4 (meaning that $K\mu_K^i\to0$ for $i>3$) thus population of trait 4 vanishes at time 0. However, Proposition \ref{CorEquilibria} implies that the rescaled macroscopic population jumps from trait $0$ to 
trait $4$ in a time 
$$s_1=\frac{-1+4/\alpha}{f_{30}}+\frac1{f_{40}}$$ 
on the $\log K$ time scale.
More precisely, the exponents are drawn on Figure \ref{fig-large-jumps}. Condition \eqref{large-jumps-cond} ensures that trait 4 fixates before trait 3.

It is easy to generalize this example to construct jumps to any distance $L$ larger than $\alpha$, by taking larger and larger fitnesses after the negative fitness region. The condition implying emergence of  trait $L$ is then a little more technical to write, since one has to compute the time for the piecewise affine function $\beta_L(t)$ (with multiple slope-breaks) to reach 1 before the other traits. Example \ref{example-large-jumps} constitutes the simplest non-trivial example of this phenomenon. Example \ref{example-large-jumps-5sites} is a further case where a more distant trait 
fixates, and two intermediate times $t_{4,1}$ and $t_{5,1}$ occur (recall the definition in \eqref{deftwk}).

\begin{example} \label{example-large-jumps-5sites}
	Let us consider the graph $\GG$ depicted on Figure \ref{fig-large-jumps-5sites}, where $V=\{0,1,2,3,4,5\}$ and $E=\{[0,1],[1,2],[2,3],[3,4],[4,5]\}$. 
	Let $\alpha\in(3,4)$, an initial condition given by $(\bar{n}(0),0,\ldots,0)$ and a fitness landscape given by
	\begin{align}
	&1,2<0,\quad 3,4,5>0,\quad 0,1,2,3,4<5 \nonumber \\
	&f_{3,0}<f_{4,0}<f_{5,0} \quad\text{and}\quad
		\frac{-1+4/\alpha}{f_{3,0}}+
		\frac{-4/\alpha+5/\alpha}{f_{4,0}}+
		\frac{1}{f_{5,0}}
		<\frac{3/\alpha}{f_{3,0}} .\label{large-jumps-cond-5sites}
	\end{align}
\end{example}

In this case, the cut-off is in between traits 3 and 4 (meaning that $K\mu_K^i\to0$ for $i>3$) thus population of trait 4 and 5 vanishes at time 0. However, 
Proposition \ref{CorEquilibria} implies that {the rescaled macroscopic population jumps from trait $0$ to trait $5$ in a time 
$$s_1=\frac{-1+4/\alpha}{f_{3,0}}+\frac{-4/\alpha+5/\alpha}{f_{4,0}}+\frac{1}{f_{5,0}}$$
on the $\log K$ time scale.
More precisely, the exponents are drawn on Figure \ref{fig-large-jumps-5sites}. Condition \eqref{large-jumps-cond-5sites} ensures that trait 5 fixates before traits 3 and 4. The first inequality is not needed but allows to simplify the second one. The dotted lines in the figures allow to construct the points where some exponents become positive.}

\begin{figure}
	\centering
	\includegraphics[width=.6\textwidth]{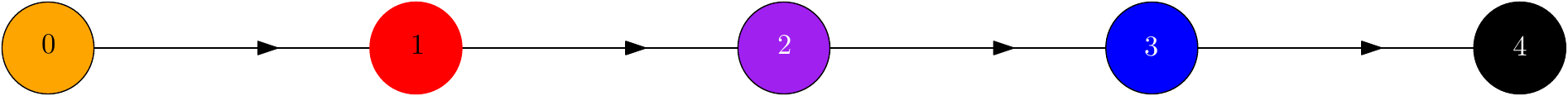}\\
	\includegraphics[width=.8\textwidth]{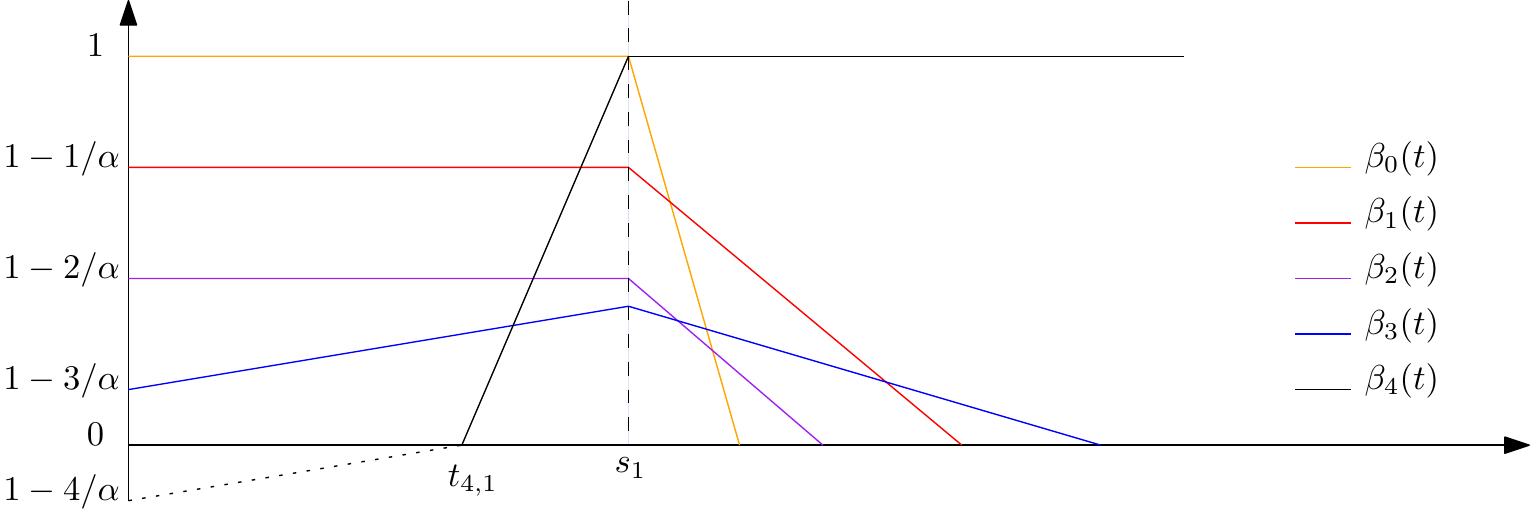}
	\caption
	{Graph $\GG$ and exponents $\beta(t)$ for Example \ref{example-large-jumps}.}
	\label{fig-large-jumps}
\end{figure}

\begin{figure}
	\centering
	\includegraphics[width=.8\textwidth,height=.6cm]{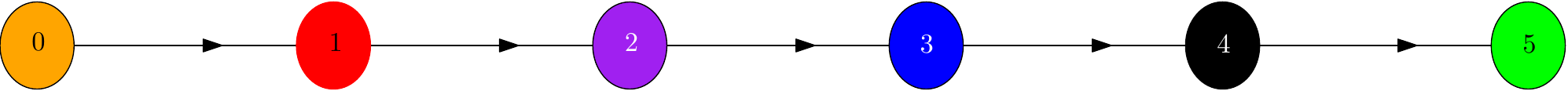}\\
	\includegraphics[width=.8\textwidth]{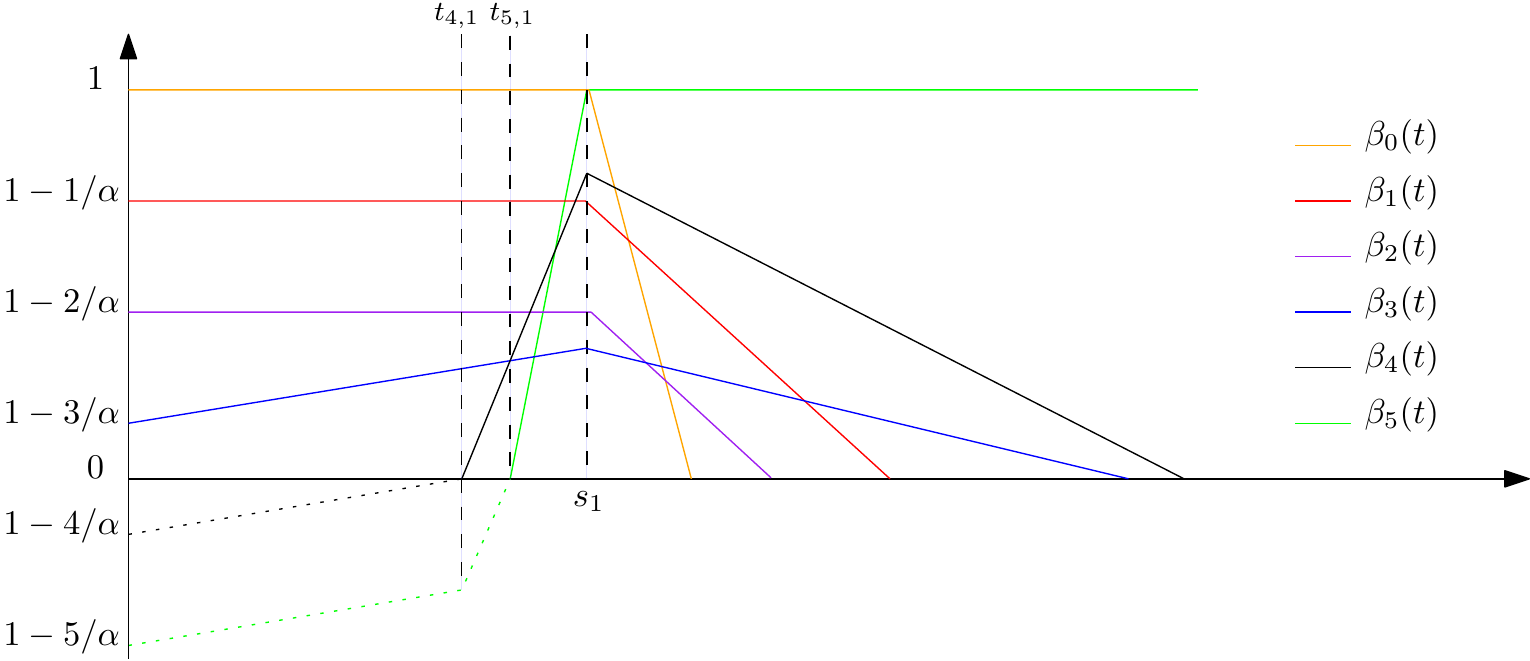}
	\caption
	{Graph $\GG$ and exponents $\beta(t)$ for Example \ref{example-large-jumps-5sites}.}
	\label{fig-large-jumps-5sites}
\end{figure}

\subsection{Effective random walk across fitness valleys}
\subsubsection{2 effective sites}

\begin{example}
	\label{example-rw2}
	Let us consider the graph $\GG$ depicted on Figure \ref{fig-rw}, where $V=\{0,1a,1b,i\}$ and $E=\{[0,i],[i,1a],[1a,1b],[i,1b]\}$. 
	We suppose that whenever there are several outgoing edges from a vertex $v$, the mutation kernel is uniform among the nearest neighboring vertices.
	Let $\alpha\in(0,1)$, an initial condition given by $(\bar{n}(0),0,\ldots,0)$ and a fitness landscape given by
	\begin{align*}
	&{
	1a\gg 0 \gg	1b\gg 1a }\\
	&{i<0,\quad i<1a,\quad i<1b. }
	\end{align*}
\end{example}

In this case, according to  \cite{BovCoqSma2018}, the time to cross the fitness valley is of order $O(1/K\mu_K^2)=O(K^{-1+2/\alpha})\gg\log K$, thus the first mutant of type $1a$ will appear on this time scale, and will invade with positive probability. Then, in a time of order $O(\log K)$, type $1b$ fixates, and one has to wait again a 
time of order $O(K^{-1+2/\alpha})$ until the appearance of the next mutant of type 0. Thus, on the time scale $O(K^{-1+2/\alpha})$, the population process converges to a jump process between the two states 0 and $1b$ with positive jump rates although the fitness $f_{1b,0}$ is negative. More precisely, following \cite{BovCoqSma2018} we define
\begin{align*}
\lambda(\rho):&=\sum_{k=1}^\infty\frac{(2k)!}{(k-1)!(k+1)!}\rho^k(1-\rho)^{k+1}\\
\rho_{i,j}:&=\frac{b_i}{b_i+d_i+c_{ij}\bar n_j}
\end{align*}

\begin{theorem} As $K\to\infty$, the following convergence holds
	$$(N^K_0,N^K_{1b})(t K^{-1+2/\alpha})\Rightarrow \bar n_{X_t}\delta_{X_t}$$
	{for finite dimensional distributions,}
	where $X_t$ is a continuous time Markov chain on $\{0,1b\}$ with transition rates
	\begin{align*}
	r_{0\to1b}&=\frac{\bar n _0 b_0}{|f_{i0}|}\lambda(\rho_{i,0})\frac{f_{1a,0}}{b_{1a}}\\
	r_{1b\to0}&=\frac{\bar n _{1b} b_{1b}}{|f_{i,1b}|}\lambda(\rho_{i,1b})\frac{f_{0,1b}}{b_{0}}.
	\end{align*}
\end{theorem}

\subsubsection{3 effective sites}

\begin{example}
	\label{example-rw3}
	Let us consider the graph $\GG$ depicted on Figure \ref{fig-rw}.
	We suppose that whenever there are several outgoing edges from a vertex $v$, the mutation kernel is uniform among the nearest neighboring vertices.
	Let $\alpha\in(0,1)$, an initial condition given by $(\bar{n}(0),0,\ldots,0)$ and a fitness landscape given by

		\begin{align*}
		\begin{array}{ccc}
		1a\gg0		&2a\gg0			&2a\gg1b\\
		1b\gg 1a	&2b\gg2a 		&1a\gg2b\\
		0\gg1b		&0\gg2b			&1b\gg2b\\
		i<0 		&j<0			&k<1a,\,k<1b\\
		i<1a,\,i<1b &j<2a,\,j<2b	&k<2a,\,k<2b
		\end{array}
		\end{align*}
\end{example}

Thus, following \cite{BovCoqSma2018}, on the time scale $O(K^{-1+2/\alpha})$, the population process converges to a jump process between the three states $\{0,1b,2b\}$ with positive jump rates. More precisely,

\begin{theorem}As $K\to\infty$, the following convergence holds
	$$(N^K_0,N^K_{1b},N^K_{2b})(t K^{-1+2/\alpha})\Rightarrow \bar n_{X_t}\delta_{X_t}$$
	{for finite dimensional distributions,} where $X_t$ is a continuous time Markov chain on $\{0,1b,2b\}$ with transition rates:

$$	\begin{array}{ll}
	r_{0\to1b}=\frac{\bar n _0 b_0/2}{|f_{i0}|}\lambda(\rho_{i,0})\frac{f_{1a,0}}{b_{1a}},&
	r_{1b\to0}=\frac{\bar n _{1b} b_{1b}/2}{|f_{i,1b}|}\lambda(\rho_{i,1b})\frac{f_{0,1b}}{b_{0}}\\
	r_{0\to2b}=\frac{\bar n _0 b_0/2}{|f_{j0}|}\lambda(\rho_{j,0})\frac{f_{2a,0}}{b_{2a}},&
	r_{2b\to0}=\frac{\bar n _{2b} b_{2b}/2}{|f_{j,2b}|}\lambda(\rho_{j,2b})\frac{f_{0,2b}}{b_{0}}\\
	r_{1b\to2b}=\frac{\bar n _{1b} b_{1b}/2}{|f_{k,1b}|}\lambda(\rho_{k,1b})\frac{f_{2a,1b}}{b_{2a}}, &
	r_{2b\to1b}=\frac{\bar n _{2b} b_{2b}/2}{|f_{k,2b}|}\lambda(\rho_{k,2b})\frac{f_{1a,2b}}{b_{1a}}.
	\end{array}$$
\end{theorem}

 \begin{figure}
	\centering
	\includegraphics[width=.3\textwidth]{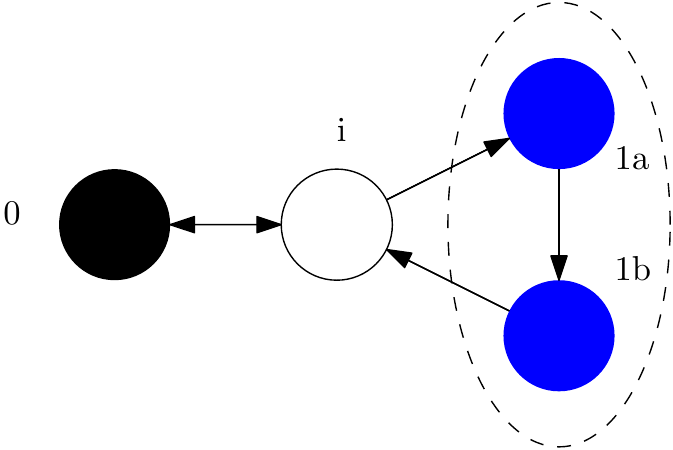} \hspace{1cm}
	\includegraphics[width=.3\textwidth]{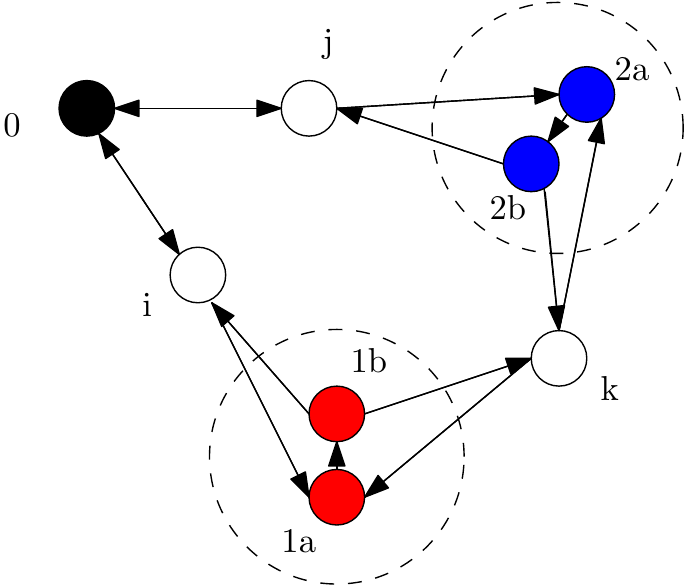}
	\caption
	{Graph $\GG$ of Examples \ref{example-rw2} and \ref{example-rw3}.}
	\label{fig-rw}
\end{figure}


\section{Proof of Theorem \ref{ThmConv} and Proposition \ref{CorEquilibria}}
\label{sec-proof}

This section is dedicated to the proofs of our main results.
As they are technical and involve many stopping times, we begin with a
rough outline of the strategy of the proof.

Throughout the proof, we define several stopping times to divide the times between invasions into sub-steps. Heuristically they correspond to the following events:
\begin{itemize}
\item $\sigma^K_k$, the time when the $k^\text{th}$ invasion has taken place and a new equilibrium is reached.
\item $\theta^K_{k,m,C}$, the first time after $\sigma^K_{k-1}$ when either the macroscopic traits stray too far from their equilibrium or at least one 
of the (formerly) microscopic traits becomes macroscopic (recall Definition \ref{def_not})
\item $s^K_k$, the first time after $\sigma^K_{k-1}$ when a microscopic trait becomes almost macroscopic, i.e. reaches an order of $K^{1-\eps_k}$.
\item $t^K_{w,k}$, the first time after $\sigma^K_{k-1}$ when trait $w$ has a positive population size. ($t^K_{w,k}=\sigma^K_{k-1}$ for all traits that are alive at this time.)
\end{itemize}
As in Proposition \ref{CorEquilibria}, $(\tau^K_\ell,\ell\geq 0)$ is the collection of both $(s^K_k,k\geq0)$ and $(t^K_{w,k},k\geq0,w\in V)$. Figure \ref{fig-stopping-times} visualises the different stopping times for the case of one macroscopic and two microscopic traits.

\begin{figure}[ht]
	\centering
	\includegraphics[width=.7\textwidth,height=6.5cm]{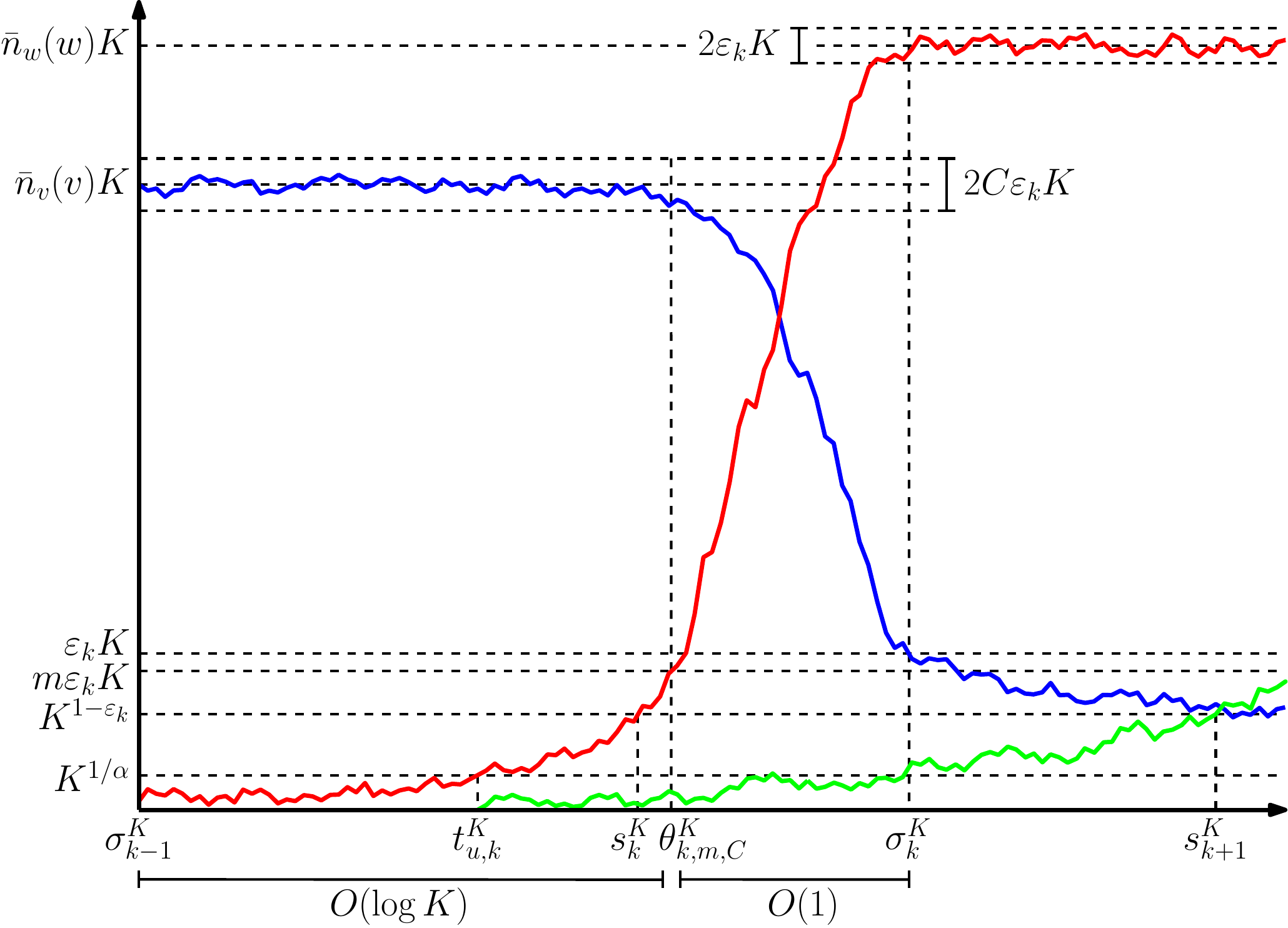}
	\caption
	{Schematic evolution of macroscopic trait $v$ (blue) and microscopic traits $w=l^K_k$ (red) and $u$ (green), where $d(w,u)=1$, during the $k^\text{th}$ invasion step.}
	\label{fig-stopping-times}
\end{figure}

The proof consists of five parts:
\begin{enumerate}
\item In the longest and most involved part of the proof, we study the growth dynamics of the different subpopulations in the time interval $[\tau^K_{\ell-1}\land T\land\theta^K_{k,m,C},\tau^K_\ell\land T\land\theta^K_{k,m,C}]$, 
making use of several results from \cite{htrans19}, which are restated in the Appendix, and generalised when needed. 
Similar to \cite{BovKra2018}, we prove lower and upper bounds for $\beta^K_w(t)$ via an induction, successively taking into account incoming mutants originating 
from traits of increasing distance to $w$. We prove that $\beta^K_w(t)$ follows the characterisation of $\beta_w(t)$ in Theorem \ref{ThmConv} up to an error 
of order $\eps_k$ for large $K$.
\item We construct the sets $M^K_\ell$ and calculate the value of $\tau^K_\ell-\tau^K_{\ell-1}$, proving part (ii) of Proposition \ref{CorEquilibria}.
\item We prove that $s^K_k$ and $\theta^K_{k,m,C}$ are equal up to an error $\eta_k$ that goes to zero as $\eps_k\to0$ and conclude that $s^K_k$ 
converges to $s_k$ when $K \to \infty$.
\item We prove that the stopping time $\theta^K_{k,m,C}$ is triggered by a (formerly) microscopic trait reaching order $K$, and not by the macroscopic traits deviating from their equilibrium.
\item Knowing that we have non-vanishing population sizes at $\theta^K_{k,m,C}$, we finally consider the Lotka-Volterra phase involving $\v_{k-1}$ and the trait $l^K_k$ that has newly reached order $K$, proving that the initial conditions for the next step, characterised in the definition of $\sigma^K_k$, are satisfied after a time of order 1. This concludes the proof of Theorem \ref{ThmConv}.
\end{enumerate}
Since part (i) of Proposition \ref{CorEquilibria} is a direct corollary of Theorem \ref{ThmConv}, this concludes the proofs of both results.\\

Recall the definitions provided in Theorem \ref{ThmConv} and Proposition \ref{CorEquilibria}, and for a given set 
$\v \subset V$, introduce $\tilde{\v}$ the support of the mutation free Lotka-Volterra equilibrium associated to $\v$, 
that is to say
$$ w \in \tilde{\v} \Leftrightarrow \bar{n}_w(\v)>0 .$$

Similarly as in \cite{htrans19}, 
the strategy of the proof consists in performing an induction on successive phases 
$k$, during which the population sizes of the set of traits $\tilde{\v}_k$ are close to their equilibrium value 
and the population sizes of the set of traits $V \setminus \tilde{\v}_k$ are small with respect to $K$.
To be more precise, we will introduce a sequence of stopping times $(\sigma_k^K \log K, k \in \N)$ 
(see definition in \eqref{def_sigma_k})
satisfying the following conditions, as soon as $s_k<T$:
\begin{assumption} \label{ass_sigma_k}\
\begin{enumerate}
 \item $\sigma_k^K \to s_k$ in probability when $K$ goes to infinity
 \item For any $0<\eps_k< 1 \wedge \inf_{w \in \tilde{\v}_k} \bar{n}_w(\v_k)$, with high probability
 \begin{enumerate}
 \item For every $w \in \tilde{\v}_k$,
 $$ \frac{N^K_w (\sigma_k^K \log K)}{K} \in \left[ \bar{n}_w(\v_k)-\eps_k,\bar{n}_w(\v_k)+\eps_k \right] $$
 \item For every $w \in \v_k \setminus \tilde{\v}_k$,
 $$ K^{1-\eps_k} \leq N^K_w (\sigma_k^K \log K) \leq \eps_k K. $$
 \item There exists $\bar{c}_k<\infty$ such that for every $w \notin \v_{k}$, either $N^K_w (\sigma_k^K \log K)=0$ if 
 $\beta_w(s_k)=0$ or
 $$ 0 < \beta_w(s_k)-\bar{c}_k \eps_k < \frac{\log \left(1+ N^K_w (\sigma_k^K \log K)\right)}{\log K}=\beta^K_w(\sigma^K_k)< \beta_w(s_k)+\bar{c}_k\eps_k<1 .$$
\end{enumerate}
 \end{enumerate} 
\end{assumption}

To be more precise, for $k \geq 1$, the time interval $[\sigma_{k-1}^K \log K,\sigma_k^K \log K]$ will be divided into two parts: 
\begin{itemize}
\item a 'stochastic phase' $[\sigma_{k-1}^K \log K,\theta_{k,m,C}^K \log K]$ needed for the trait 
$$ l^K_{k}:= \v^K_k \setminus \v^K_{k-1} $$
to reach a size of order $K$, 
\item a 'deterministic phase' $[\theta_{k,m,C}^K \log K,\sigma_{k}^K \log K]$ needed for
the mutation free Lotka-Volterra system associated to $\tilde{\v}^K_{k-1} \cup l^K_k$ to reach a neighbourhood of its equilibrium.
\end{itemize}

\subsection*{Initialisation of the induction.}
\noindent $\bullet$ \underline{$\sigma_0^K$:} 
By assumption, 
$$ \beta^K_w(0)\overset{K\to\infty}{\rightarrow}\left(1-\frac{d(\v_0,w)}{\alpha}\right)_+. $$
Let us choose a small $\eps_0>0$.
Then from point $(ii)$ of Lemma \ref{lemmaC1}, there exists {a deterministic} $T(\eps_0)<\infty$ such that
$$ \lim_{K \to \infty} \P \left( \| N^K(T(\eps_0))/K- \bar{n}(\v_0) \|_\infty \leq \eps_0  \right)=1. $$
Define $\sigma_0^K:=T(\eps_0)/\log K$.
We can check that $\sigma_0^K$ is a stopping time converging in probability to $s_0=0$ and satisfying 
Assumption \ref{ass_sigma_k}.
Moreover we know that the processes $\beta_w^K, w \in V$, vary on a time scale of order $\log K$ (see \cite{Cha06,htrans19} for instance). In particular, 
they do not vary during the time $T(\eps_0)$ in the large $K$ limit.
This entails that $\sigma_0^K$ satisfies Assumption \ref{ass_sigma_k}.
\\

\noindent $\bullet$ \underline{$\sigma_k^K, k \geq 1$:} 

Assume that $s_{k-1}<T_0$ and that $\sigma_{k-1}^K \log K$ 
is a stopping time satisfying Assumption \ref{ass_sigma_k}. We will now construct $\sigma_{k}^K$.

\subsection{Definitions and first properties}

Let us introduce a small $\eps_{k}>0$ as well as a stopping time $\theta_{k,m,C}^K \log K $ via 
\begin{equation} \label{def_theta_k} \theta_{k,m,C}^K:= \inf \left\{ t \geq \sigma_{k-1}^K, \exists w \in \v_{k-1}, 
\left| \frac{N_w^K(t \log K)}{K}-\bar{n}_w(\v_{k-1}) \right| \geq C \eps_k \text{ or }
\sum_{w \notin \v_{k-1}}N_w^K(t \log K) \geq m \eps_k K
\right\}. \end{equation}
The conditions satisfied by $m>0$ and $C>0$ will be precised later on. $m$ is typically small, see \eqref{stay_close}.
The conditions satisfied by $C$ will be specified in Section \ref{pageC}.

We will now finely study the population dynamics on the time interval 
$[\sigma_{k-1}^K \log K ,(\theta_{k,m,C}^K \wedge T) \log K]$. To this aim, 
we will couple the subpopulations of individuals with a given trait with branching processes 
with immigration
and use results on these processes derived in \cite{htrans19} and recalled (and generalized when needed) in the Appendices.
The main difficulty of this step comes from the fact that as we allow for any finite graph of mutations, 
the immigration rate for a particular subpopulation may vary a lot on the time interval 
$[\sigma_{k-1}^K \log K,(\theta_{k,m,C}^K\wedge T) \log K]$. 
This is why we introduced in Proposition \ref{CorEquilibria} the sequence of times $(\tau_\ell, \ell \in \N)$, which corresponds to 
the times when mutants of a new type arise or a formerly microscopic trait becomes of order K.

Notice that although we make extensive use of the techniques and results developed in \cite{htrans19}, the authors of this paper considered a specific graph 
embedded in $\Z$, and their proof structure, in particular inductions, relies on their graph structure.
The current inductions are more involved and more in the proof spirit of \cite{BovKra2018}.
\\

To begin with, let us recall the rates of the different events for the population $N_w^K$, with $w \in V$, at time 
$t$:

\begin{itemize}
 \item[$\bullet$] Reproductions without mutation:
\begin{equation}\label{birth_rate_w}
 \mathfrak{b}_w(t):= b_w (1-K^{-1/\alpha})N_w^K(t)
\end{equation}

 \item[$\bullet$] Death:
 \begin{equation}\label{death_rate_w}
 \mathfrak{d}_w(t):= \left(d_w + \sum_{x \in V}\frac{c_{w,x}}{K}N_x^K(t) \right)N_w^K(t)
\end{equation}
 \item[$\bullet$] Reproductions with mutations towards the trait $w$:
 \begin{equation}\label{birthmut_rate_w}
 \mathfrak{bm}_w(t):= K^{-1/\alpha}\sum_{x \in V, d(x,w)=1} b_x 
 m(x,w) N_x^K(t).
\end{equation}
\end{itemize}

Notice that for $K$ large enough, as $\sigma_{k-1}^K$ satisfies Assumption \ref{ass_sigma_k} and 
by definition of $\theta_{k,m,C}^K$, on the time interval $[\sigma_{k-1}^K \log K,(\theta_{k,m,C}^K \wedge T) \log K]$, we have

	\begin{align} 
		\label{bound_b}  b(w,k,{-})N_w^K(t) &\leq
		\mathfrak{b}_w(t)\leq   b(w,k,{+})N_w^K(t),\\
		\label{bound_d}  d(w,\tilde{\v}_{k-1},k,+)N_w^K(t) &\leq
		\mathfrak{d}_w(t)\leq   d(w,\tilde{\v}_{k-1},k,-)N_w^K(t)\\
		\label{bound_f}f(w,\tilde{\v}_{k-1},k,-)&\leq f_{w,\tilde{\v}_{k-1}}\leq f(w,\tilde{\v}_{k-1},k,+),
		\end{align}
where we have introduced the following notations, for any $w \in V$ and $* \in \{-,+\}$,
\begin{align}
 b(w,k,-)&:= (1-\eps_k)b_w, \quad b(w,k,+):= b_w, \nonumber \\
 d(w,\tilde{\v}_{k-1},k,-)&:= d_w + \sum_{x \in \tilde{\v}_{k-1}}c_{w,x}\bar{n}_x(\v_{k-1}) 
+ \left(\sum_{x \in V}c_{w,x}\right)(m+C)\eps_k, \nonumber \\
 d(w,\tilde{\v}_{k-1},k,+)&:= d_w + \sum_{x \in \tilde{\v}_{k-1}}c_{w,x}\bar{n}_x(\v_{k-1}) - \left(\sum_{x \in V}c_{w,x}\right)C\eps_k, \nonumber \\
	f(w,\tilde{\v}_{k-1},k,*)&:= b(w,k,*)-d(w,\tilde{\v}_{k-1},k,*). \label{def_f}  
\end{align}

Hence the rate of reproduction without mutation, as well as the death rate do not vary significantly during the time 
interval $[\sigma_{k-1}^K \log K,(\theta_{k,m,C}^K\wedge T) \log K]$. The difficulty comes from the rate of mutations towards a given trait, 
which depends on the population sizes of its neighbours in the graph $\mathcal{G}$, which themselves depend on the population sizes
of their neighbours and so on.\\

Let us introduce the times $\tau_\ell^K$ and the sets $M^K_\ell$, which correspond respectively to the 
times of invasion or appearance of new mutants (and will be the time steps of the algorithm to be described shortly later)
and to the sets of living traits {in the time interval $(\tau_\ell^K,\tau_{\ell+1}^K]$}.
To be more precise, 
\begin{definition}
	Let $ s_k^K:= \inf \{ t \geq \sigma_{k-1}^K: \exists w \in V \setminus \v^K_{k-1}, \beta_w^K(t)> 1- \eps_k \}, $
	and
	\begin{align*}
	t^K_{w,k}:=\begin{cases}\inf\{t\geq \sigma^K_{k-1}:\exists\ u\in V: d(u,w)=1, \beta^K_u(t)=\frac{1}{\alpha}\}
	&\text{if }\beta_w^K(\sigma^K_{k-1})=0\\\sigma^K_{k-1}&\text{else},\end{cases}
	\end{align*}
	 {The sequences $(\tau^K_\ell, \ell\geq0)$ and $(M^K_\ell,\ell\geq0)$} are defined as follows:\\
	$\tau^K_0=\sigma_0^K$ and, for $\sigma^K_{k-1}\leq \tau_{\ell-1}^K<s_k^K$,
	$$ \tau_\ell^K=s_k^K\land\min \ \{t^K_{w,k}:w\in V, t^K_{w,k}>\tau^K_{\ell-1}\}, $$
that is to say the minimum between the time when a previously microscopic population becomes (almost) macroscopic, and the time of appearance 
of a new mutant.
From the definition of the sequence $(\tau^K_\ell, \ell\geq0)$ we can now define the sequence of sets of living traits $(M^K_\ell,\ell\geq0)$ via
	\begin{align*}
	M_\ell^K&=\{w\in V:\beta^K_w(\tau^K_\ell)>0\text{ or }\tau^K_\ell=t^K_{w,k}\}\notag\\
	&=\big(M^K_{\ell-1}\backslash\{w\in V:\beta^K_w(\tau^K_\ell)=0\}\big)\cup\{w\in V:\tau_\ell^K=t^K_{w,k}\}.
	\end{align*}
\end{definition}

\subsection{Dynamics of the process on $[\tau^K_{\ell-1}\log K,\tau^K_{\ell}\log K]$}
\emph{}\\

We will first
prove that 
there exists a finite and positive constant $C_\ell$ such that with high probability, 
for every $w\in M^K_{\ell-1}$ and $t\in 
[\tau^K_{\ell-1} \wedge T \wedge \theta_{k,m,C}^K,\tau^K_{\ell} \wedge T \wedge \theta_{k,m,C}^K]$, 
\begin{multline} \label{inductiontoprove}
\max_{\substack{u\in M^K_{\ell-1}}}\left[\beta^K_u(\tau^K_{\ell-1})-\frac{d(u,w)}{\alpha}+(t-\tau^K_{\ell-1})
f(u,\tilde{\v}_{k-1},k,-)\right]_+
\\
\leq \beta^K_w(t) \leq \\
\max_{\substack{u\in M^K_{\ell-1}}}\left[\beta^K_u(\tau^K_{\ell-1})-\frac{d(u,w)}{\alpha}+ C_\ell \eps_k+(t-\tau^K_{\ell-1})
f(u,\tilde{\v}_{k-1},k,+)\right]_+\end{multline}

Let us thus take $t$ in $[\tau^K_{\ell-1} \wedge T \wedge \theta_{k,m,C}^K,\tau^K_{\ell} \wedge T \wedge \theta_{k,m,C}^K]$.
To obtain the lower bound in \eqref{inductiontoprove}, we show by induction that, for any $n \geq 0$ and with high probability,

\begin{equation}\label{LB-beta}
\beta^K_w(t) \geq \max_{\substack{u\in M^K_{\ell-1}:\\d(u,w)\leq n}}
\left[\beta^K_u(\tau^K_{\ell-1})-\frac{d(u,w)}{\alpha}+(t-\tau^K_{\ell-1})
f(u,\tilde{\v}_{k-1},k,-)\right]_+
\end{equation}


\noindent \textbf{Induction lower bound:}
\noindent $\bullet$ \underline{$n=0$:} let $w\in M^K_{\ell-1}$. From \eqref{bound_b} and \eqref{bound_d}, 
we see that we can couple $N^K_w$ with a process 
$Z^K$ with law $BP_K\left(b(w,k,-),d(w,\tilde{\v}_{k-1},k,{-}),\beta^K_w(\tau^K_{\ell-1})\right)$ 
(see the definition of $BP_K$ in Section \ref{def_BPK})
in such a way that
$$
N^K_w(t\log K)\geq Z^K((t-\tau^K_{\ell-1})\log K).
$$
Hence, from Corollary \ref{thm-line-break}, we obtain that {with high probability}, 
\begin{align*}
\beta^K_w(t)\geq \left[\beta^K_w(\tau^K_{\ell-1})+(t-\tau^K_{\ell-1}){f(w,\tilde{\v}_{k-1},k,-)}\right]_+.
\end{align*}

\begin{remark}\label{rem-[eps,T]}
Notice that the application of Lemma \ref{lemmaA1} (which has been derived in \cite{htrans19}) would require $\beta^K_w(\tau^K_{\ell-1})>0$ and that this condition {may not be} satisfied 
for one of the $w \in \tilde{\v}_{k-1}$ (the trait which becomes macroscopic at time $\tau^K_{\ell-1} \log K$). However, the population of individuals $w$ 
grows exponentially due to the mutations coming from another trait and there exists a finite $c$ such that, for small $\delta>0$, $N_w^K((\tau^K_{\ell-1}+\delta) \log K) \geq K^{c\delta}$.
We could thus apply Lemma \ref{lemmaA1} at this time, and later on let $\delta$ go to $0$ to get the result. 
This is in words the statement of Corollary \ref{thm-line-break}.
\end{remark}

\noindent $\bullet$ \underline{$n\to n+1$:} 
Let $w,u',u\in M^K_{\ell-1}$ such that $d(u',w)=1$ and $d(u,u')\leq n$.
From now on, we will use the notation $BPI_K$, which is defined in Section \ref{sect_BPI}.
From \eqref{bound_b}, \eqref{bound_d}, and \eqref{bound_f}, 
by looking only at the immigration coming from $u'$, we see that we can couple $N^K_w$ with a process 
$Z^K$ with law 
$$
BPI_K\left(b(w,k,-),d(w,\tilde{\v}_{k-1},k,-),f(u',\tilde{\v}_{k-1},k,-),\beta^K_{u'}(\tau^K_{\ell-1})-
\frac{1}{\alpha},\beta^K_w(\tau^K_{\ell-1})\right)
$$
in such a way that
$$
N^K_w(t\log K)\geq Z^K((t-\tau^K_{\ell-1})\log K).
$$
By the induction hypothesis, with high probability, 
\begin{align}\label{HR}
\beta^K_{u'}(t)\geq\left[\beta^K_u(\tau^K_{\ell-1})-\frac{d(u,u')}{\alpha}+(t-\tau^K_{\ell-1}){f(u,\tilde{\v}_{k-1},k,{-})}\right]_+,
\end{align}
which implies that we can couple $Z^K$ with a process $Y^K$ with law
$$BPI_K\left(b(w,k,-),d(w,\tilde{\v}_{k-1},k,-),f(u,\tilde{\v}_{k-1},k,-),\beta^K_u(\tau^K_{\ell-1})-
\frac{d(u,u')+1}{\alpha},\beta^K_w(\tau^K_{\ell-1})\right)$$
in such a way that
$$
Z^K((t-\tau^K_{\ell-1})\log K)\geq Y^K((t-\tau^K_{\ell-1})\log K).
$$

Hence, from Corollary \ref{thm-line-break}, even if we have to work in a time interval $[\tau_{\ell-1}^K+\delta,T]$,
for a small positive $\delta$,
in the spirit of Remark \ref{rem-[eps,T]}, 
as $w \in M^K_{\ell-1}$ 
we obtain that with high probability,

\begin{align*}
\beta^K_{w}(t)&\geq 
\left[\beta^K_w(\tau^K_{\ell-1})\lor \left( \beta^K_u(\tau^K_{\ell-1})-\frac{d(u,u')+1}{\alpha} \right)
+(t-\tau^K_{\ell-1}){f(w,\tilde{\v}_{k-1},k,-)}\right]_+\\
&\hspace{6cm}\lor\left[\beta^K_u(\tau^K_{\ell-1})-\frac{d(u,u')+1}{\alpha}+(t-\tau^K_{\ell-1}){f(u,\tilde{\v}_{k-1},k,-)}\right]_+\\
&\geq \left[\beta^K_w(\tau^K_{\ell-1})
+(t-\tau^K_{\ell-1}){f(w,\tilde{\v}_{k-1},k,-)}\right]_+\\
&\hspace{6cm}\lor\left[\beta^K_u(\tau^K_{\ell-1})-\frac{d(u,u')+1}{\alpha}+(t-\tau^K_{\ell-1}){f(u,\tilde{\v}_{k-1},k,-)}\right]_+.
\end{align*}

As this is true for any $u'$ such that $d(u',w)=1$ and as the above bound is a decreasing function of $d(u,u')$, by taking the supremum over such $u'$ we obtain 
\begin{multline*}
\beta^K_{w}(t)\geq 
\left[\beta^K_w(\tau^K_{\ell-1})+(t-\tau^K_{\ell-1}){f(w,\tilde{\v}_{k-1},k,-)}\right]_+\\
\lor\left[\beta^K_u(\tau^K_{\ell-1})-\frac{d(u,w)}{\alpha}+(t-\tau^K_{\ell-1}){f(u,\tilde{\v}_{k-1},k,-)}\right]_+.
\end{multline*}

Thus, with high probability,

\begin{align*}
\beta^K_w(t)\geq\max_{\substack{u\in M^K_{\ell-1}:\\d(u,w)\leq n+1}}\left[\beta^K_u(\tau^K_{\ell-1})-\frac{d(u,w)}{\alpha}+
(t-\tau^K_{\ell-1}){f(u,\tilde{\v}_{k-1},k,-)}\right]_+,
\end{align*}

which ends the induction for the lower bound.\\

Let us now proceed to the induction for the upper bound.
We again take $t$ in $[\tau^K_{\ell-1} \wedge T \wedge \theta_{k,m,C}^K,\tau^K_{\ell} \wedge T \wedge \theta_{k,m,C}^K]$
and we will show that {for any $n \in \N$ there exists a finite constant $C_{n,\ell}$ such that} {with high probability},
\begin{multline*}\beta^K_w(t) \leq 
\max_{\substack{u\in M^K_{\ell-1}:\\d(u,w)\leq n}}\left[{ \left(\beta^K_u(\tau^K_{\ell-1})-\frac{d(u,w)}{\alpha} {+ C_{n,\ell} \eps_k} \right)
\lor\left(1-\frac{n+1}{\alpha}+(n+2)\eps_k\right)}+(t-\tau^K_{\ell-1}){f(u,\tilde{v}_{k-1},k,+)}\right]_+\\
\lor\left(1-\frac{n+1}{\alpha}+(n+2)\eps_k\right).
\end{multline*}

 Notice that, since $\eps_k$ can be chosen small enough such that $1-(n+1)/\alpha+(n+2)\eps_k<0$, 
 for all $n>\lfloor\alpha\rfloor$, all terms with $d(u,w)>\lfloor\alpha\rfloor$ are negative and this equation is equivalent 
 to the upper bound in \eqref{inductiontoprove}.\\

\noindent \textbf{Induction upper bound:}

Throughout the induction for the upper bound, we will several times make use of the fact that we can approximate the total 
immigration to one trait, which is the sum of the mutants coming from its neighbours, from above by the number of neighbours times the largest incoming mutation. 
More precisely, if $I_w$ is the number of incoming neighbours of $w$,
\begin{align}\label{sumapprox}
\sum_{\substack{u\in V:\\d(u,w)=1}}K^{\beta^K_u(t)}\leq I_w \max_{\substack{u\in V:\\d(u,w)=1}}K^{\beta^K_u(t)}=
\max_{\substack{u\in V:\\d(u,w)=1}}K^{(\log I_w/\log K)+\beta^K_u(t)}.
\end{align}
Since the trait space is finite, for $K$ large enough, we can assume that $\max_{w\in V}\log I_w/\log K\leq \eps_k$.

\noindent $\bullet$ \underline{$n=0$:} We observe that for $K$ large enough $\beta^K_u(t)\leq1+\eps_k$ for every $u\in M^K_{\ell-1}$ 
such that $d(u,w)=1$ (see Corollary \ref{betabound}).

From \eqref{bound_b}, \eqref{bound_d}, and \eqref{bound_f}, 
we see that we can couple $N^K_w$ with a process 
$Z^K$ with law $BPI_K\left({b(w,k,+),d(w,\tilde{\v}_{k-1},k,+)},0,1-\frac{1}{\alpha}+2\eps_k,\beta^K_w(\tau^K_{\ell-1})\right)$ 
in such a way that
$$
N^K_w(t\log K)\leq Z^K((t-\tau^K_{\ell-1})\log K).
$$
Hence from Corollary \ref{thm-line-break}, even if we have to work in a time interval 
$[\tau_{\ell-1}^K+\delta,T]$, for a small positive $\delta$, in the spirit of Remark \ref{rem-[eps,T]}, as {$w \in M^K_{\ell-1}$} 
we obtain that with high probability,
\begin{align*}
\beta^K_w(t)\leq\left[\beta^K_w(\tau^K_{\ell-1})
{\lor\left(1-\frac{1}{\alpha}+2\eps_k\right)}
+(t-\tau^K_{\ell-1}){f(w,\tilde{v}_{k-1},k,+)}\right]_+\lor\left(1-\frac{1}{\alpha}+2\eps_k\right).
\end{align*}

\noindent $\bullet$ \underline{$n\to n+1$:} For $w,u'\in M^K_{\ell-1}$ such that $d(u',w)=1$, by the induction 
hypothesis we have {the existence of a finite constant $C_{n,\ell}$ such that}, with high probability,
\begin{multline*}
\beta^K_{u'}(t)\leq\max_{\substack{u\in M^K_{\ell-1}:\\d(u,u')\leq n}}
\left[{ \left(\beta^K_u(\tau^K_{\ell-1})-\frac{d(u,u')}{\alpha}+ C_{n,\ell} \eps_k\right) \lor\left(1-\frac{n+1}{\alpha}+(n+2)\eps_k\right)}+
(t-\tau^K_{\ell-1}){f(u,\tilde{v}_{k-1},k,+)}\right]_+
\\ \lor\left(1-\frac{n+1}{\alpha}+(n+2)\eps_k\right).
\end{multline*}

From \eqref{bound_b}, \eqref{bound_d}, and \eqref{bound_f}, by looking at the maximal immigration coming from a neighbouring $u'$ and adding another $\eps_k$ in the spirit of \eqref{sumapprox},
we thus see that we can couple $N^K_w$ with multiple processes 
$Z^{K,u,u'}$ and $Z^K$ with respective laws 
\begin{multline*}BPI_K\Big({b(w,k,+),d(w,\tilde{\v}_{k-1},k,+)},{f(u,\tilde{v}_{k-1},k,+)},\\
{ \left(\beta^K_u(\tau^K_{\ell-1})-\frac{d(u,u')+1}{\alpha}+ (C_{n,\ell}+1) \eps_k\right) \lor\left(1-\frac{n+2}{\alpha}+(n+3)\eps_k\right)},\beta^K_w(\tau^K_{\ell-1})\Big)\end{multline*}
and
$$BPI_K\left({b(w,k,+),d(w,\tilde{\v}_{k-1},k,+)},0,1-\frac{n+2}{\alpha}+(n+3)\eps_k,\beta^K_w(\tau^K_{\ell-1})\right)$$
in such a way that
$$N^K_w(t\log K)\leq \max_{\substack{u'\in M^K_{\ell-1}:\\d(u',w)=1}}
\max_{\substack{u\in M^K_{\ell-1}:\\d(u,u')\leq n}}Z^{K,u,u'}((t-\tau^K_{\ell-1})\log K)\lor Z^K((t-\tau^K_{\ell-1})\log K).
$$

Hence from Corollary \ref{thm-line-break}, even if we have to work in a time interval 
$[\tau_{\ell-1}^K+\delta,T]$, for a small positive $\delta$, in the spirit of Remark \ref{rem-[eps,T]}, as $w \in M^K_{\ell-1}$
we obtain that {with high probability},
{\scriptsize \begin{align}
\beta^K_{w}(t)\leq&\max_{\substack{u'\in M^K_{\ell-1}:\\d(u',w)=1}}\max_{\substack{u\in M^K_{\ell-1}:\\d(u,u')\leq n}}
\left\{\left[\beta^K_w(\tau^K_{\ell-1}) { \lor \left(\beta^K_u(\tau^K_{\ell-1})-\frac{d(u,u')+1}{\alpha}+ (C_{n,\ell}+1) \eps_k\right) \lor\left(1-\frac{n+2}{\alpha}+(n+3)\eps_k\right)}+(t-\tau^K_{\ell-1}){f(w,\tilde{v}_{k-1},k,+)}\right]_+\notag \right. \\
& \left.\lor\left[{ \left(\beta^K_u(\tau^K_{\ell-1})-\frac{d(u,u')+1}{\alpha}+ (C_{n,\ell}+1) \eps_k\right) \lor\left(1-\frac{n+2}{\alpha}+(n+3)\eps_k\right)}+(t-\tau_{\ell-1}){f(u,\tilde{v}_{k-1},k,+)}\right]_+
\lor\left(1-\frac{n+2}{\alpha}+(n+3)\eps_k\right)\right\}  \nonumber \\
\leq&\max_{\substack{u\in M^K_{\ell-1}:\\d(u,w)\leq n+1}}
\left\{\left[\beta^K_w(\tau^K_{\ell-1}) { \lor \left(\beta^K_u(\tau^K_{\ell-1})-\frac{d(u,w)}{\alpha}+ (C_{n,\ell}+1) \eps_k\right) \lor\left(1-\frac{n+2}{\alpha}+(n+3)\eps_k\right)}+(t-\tau^K_{\ell-1}){f(w,\tilde{v}_{k-1},k,+)}\right]_+ \right. \notag\\
& \left. \lor\left[{ \left(\beta^K_u(\tau^K_{\ell-1})-\frac{d(u,w)}{\alpha}+ (C_{n,\ell}+1) \eps_k\right) \lor\left(1-\frac{n+2}{\alpha}+(n+3)\eps_k\right)}+(t-\tau_{\ell-1}){f(u,\tilde{v}_{k-1},k,+)}\right]_+
\lor\left(1-\frac{n+2}{\alpha}+(n+3)\eps_k\right)\right\}. \label{LBbeta2}
\end{align}}

In order to simplify the right hand side of the previous inequality, we will show that for any $\ell \in \N$ there exists a finite and positive constant $C_{\ell}$ such that for any $(u,w) \in V^2$, with high probability
\begin{equation} \label{ineq2_beta}  \beta^K_u(\tau^K_{\ell-1})-\frac{d(u,w)}{\alpha} \leq \beta^K_w(\tau^K_{\ell-1})+ C_{\ell} \eps_k. \end{equation}

Combining \eqref{LBbeta2} and \eqref{ineq2_beta} yields that {with high probability},
{\small\begin{align*}
\beta^K_{w}(t)\leq&\max_{\substack{u\in M^K_{\ell-1}:\\d(u,w)\leq n+1}}
\left\{\left[\left(\beta^K_w(\tau^K_{\ell-1})+ (C_{n,\ell}+1+C_\ell) \eps_k \right) \lor\left(1-\frac{n+2}{\alpha}+(n+3)\eps_k\right)+(t-\tau^K_{\ell-1}){f(w,\tilde{v}_{k-1},k,+)}\right]_+ \right. \notag\\
& \left. \lor\left[{ \left(\beta^K_u(\tau^K_{\ell-1})-\frac{d(u,w)}{\alpha}+ (C_{n,\ell}+1) \eps_k\right) \lor\left(1-\frac{n+2}{\alpha}+(n+3)\eps_k\right)}+(t-\tau_{\ell-1}){f(u,\tilde{v}_{k-1},k,+)}\right]_+\right.\notag\\
& \left.\lor\left(1-\frac{n+2}{\alpha}+(n+3)\eps_k\right)\right\}\\
\leq & \max_{\substack{u\in M^K_{\ell-1}:\\d(u,w)\leq n+1}}
 \left[{ \left(\beta^K_u(\tau^K_{\ell-1})-\frac{d(u,w)}{\alpha}+ (C_{n,\ell}+1+C_\ell)\eps_k\right) \lor\left(1-\frac{n+2}{\alpha}+(n+3)\eps_k\right)}+(t-\tau_{\ell-1}){f(u,\tilde{v}_{k-1},k,+)}\right]_+
\\& \qquad \lor\left(1-\frac{n+2}{\alpha}+(n+3)\eps_k\right),
\end{align*}}
which ends the induction for the upper bound.\\
Let us now derive inequality \eqref{ineq2_beta}. It is obtained by an induction on $\ell$. 
If $\ell=1$, by \eqref{cond-init} and the triangle inequality,
\begin{align*}
 \lim_{K \to \infty} \beta_u^K(0) - \frac{d(u,w)}{\alpha} & = \left[ 1 - \frac{d(\v_0,u)}{\alpha} \right]_+  - \frac{d(u,w)}{\alpha} \\
 & \leq \left[ 1 - \frac{d(\v_0,u)}{\alpha}- \frac{d(u,w)}{\alpha} \right]_+  
 \leq \left[ 1 - \frac{d(\v_0,w)}{\alpha} \right]_+
 =  \lim_{K \to \infty} \beta_w^K(0) .
\end{align*}
As the convergence is in probability, it means that for $K$ large enough, there exists a finite $C_{u,w}$ such that with a probability larger than $1- \eps_k$,
\begin{equation} \label{inductionl0} \beta_u^K(0) - \frac{d(u,w)}{\alpha} \leq  \beta_w^K(0) + C_{u,w}\eps_k. \end{equation}
As there are only finitely many traits, $\sup_{u,w\in V} C_{u,w} < \infty$.
Moreover, as $\eps_k$ can be chosen as small as we want and as we want to prove a convergence in probability, we may focus on the event where inequality \eqref{inductionl0} is satisfied. We will do that later on without mentioning it again for the sake of readability.

Now assume that \eqref{ineq2_beta} is true for $\ell-1 \in \N$.
Let us first prove that it still holds for $\ell$.

From the previous step on the time interval $[\tau^K_{\ell-2} \wedge T \wedge \theta_{k,m,C}^K,\tau^K_{\ell-1} \wedge T \wedge \theta_{k,m,C}^K]$,
we know that if $\tau^K_{\ell-1} \leq T \wedge \theta_{k,m,C}^K$, for any $w \in V$ and $K$ large enough,
$$ \max_{u\in M^K_{\ell-2}}\left[\beta^K_u(\tau^K_{\ell-2})-\frac{d(u,w)}{\alpha}+(\tau^K_{\ell-1}-\tau^K_{\ell-2})
f(u,\tilde{\v}_{k-2},k,-)\right]_+
\\
\leq \beta^K_w(\tau^K_{\ell-1}). $$
Now let us take $u \in V$. We also deduce from the previous step that for $K$ large enough
$$ \beta^K_u(\tau^K_{\ell-1}) \leq  \max_{u'\in M^K_{\ell-2}}
\left[\beta^K_{u'}(\tau^K_{\ell-2})-\frac{d(u',u)}{\alpha}+ C_{\ell-1} \eps_k +(\tau^K_{\ell-1}-\tau^K_{\ell-2})
f(u',\tilde{\v}_{k-2},k,+)\right]_+. $$
In particular there exists $\tilde{u} \in V$ such that $d(\tilde{u},u)\leq \lfloor \alpha \rfloor$ and for $K$ large enough
$$ \beta^K_u(\tau^K_{\ell-1}) \leq \left[\beta^K_{\tilde{u}}(\tau^K_{\ell-2})-\frac{d(\tilde{u},u)}{\alpha}+ C_{\ell-1} \eps_k+(\tau^K_{\ell-1}-\tau^K_{\ell-2})
f(\tilde{u},\tilde{\v}_{k-2},k,+)\right]_+. $$
Thus, for $K$ large enough,
\begin{align*} \beta^K_u(\tau^K_{\ell-1})- \frac{d(u,w)}{\alpha} & 
\leq \left[\beta^K_{\tilde{u}}(\tau^K_{\ell-2})-\frac{d(\tilde{u},u)}{\alpha}+ C_{\ell-1} \eps_k+(\tau^K_{\ell-1}-\tau^K_{\ell-2})
f(\tilde{u},\tilde{\v}_{k-2},k,+)\right]_+ - \frac{d(u,w)}{\alpha}\\
& \leq \left[\beta^K_{\tilde{u}}(\tau^K_{\ell-2})-\frac{d(\tilde{u},u)+d(u,w)}{\alpha}+ C_{\ell-1} \eps_k+(\tau^K_{\ell-1}-\tau^K_{\ell-2})
f(\tilde{u},\tilde{\v}_{k-2},k,+)\right]_+ \\
& \leq \left[\beta^K_{\tilde{u}}(\tau^K_{\ell-2})-\frac{d(\tilde{u},w)}{\alpha}+ C_{\ell-1} \eps_k+(\tau^K_{\ell-1}-\tau^K_{\ell-2})
f(\tilde{u},\tilde{\v}_{k-2},k,+)\right]_+\\
& \leq \left[\beta^K_{\tilde{u}}(\tau^K_{\ell-2})-\frac{d(\tilde{u},w)}{\alpha}+ C_{\ell-1} \eps_k+(\tau^K_{\ell-1}-\tau^K_{\ell-2})
f(\tilde{u},\tilde{\v}_{k-2},k,-)\right]_+ + C \eps_k \\
& \leq  \max_{\tilde{u}\in M^K_{\ell-2}} \left[\beta^K_{\tilde{u}}(\tau^K_{\ell-2})-\frac{d(\tilde{u},w)}{\alpha}
+(\tau^K_{\ell-1}-\tau^K_{\ell-2})
f(\tilde{u},\tilde{\v}_{k-2},k,-)\right]_+ + (C_{\ell-1} + C) \eps_k \\
& \leq  \beta^K_{w}(\tau^K_{\ell-1})+ (C_{\ell-1} + C) \eps_k, \end{align*}
where we used \eqref{def_f}, \eqref{LB-beta}, the bound $\tau^K_{\ell-1}-\tau^K_{\ell-2}\leq T$, and
$$ C :=  \max_{\substack{\tilde{u}\in M^K_{\ell-2}}} \left( b_{\tilde{u}}+ (m+6) \left( \sum_{x \in V} c_{\tilde{u},x} \right) \right)T.$$
This entails \eqref{ineq2_beta}.\\

To conclude the proof of \eqref{inductiontoprove}, we just need to notice that for 
$n>\lfloor \alpha \rfloor$, if $\eps_k$ is small enough, 
$$ 1-\frac{n+1}{\alpha}+(n+2)\eps_k<0. $$

As $\tau^K_{\ell}- \tau^K_{\ell-1}\leq T$, Equation \eqref{inductiontoprove} tells us that, with an error of order $\eps_k$ which is as small as we want, {with high probability}, the growth of traits $w\in M^K_{\ell-1}$ follows, for 
$t\in[\tau^K_{\ell-1} \wedge T \wedge \theta_{k,m,C}^K,\tau^K_{\ell} \wedge T \wedge \theta_{k,m,C}^K]$
\begin{align*}
\beta^K_w(t)\cong\max_{u\in M^K_{\ell-1}}\left[\beta_u^K(\tau^K_{\ell-1})-\frac{d(u,w)}{\alpha}+
(t-\tau^K_{\ell-1})f_{u,\v_{k-1}}\right]_+.
\end{align*}
To avoid repetition, we will write $\cong$ in the sequel to indicate approximations with high probability, with an error of order $\eps_k$. 

\subsection{Value of $\tau^K_\ell$ and construction of $M^K_\ell$}

Let us assume for the moment (it will be proven in Section \ref{pageC}) that the following holds with high probability:
\begin{equation} \label{thetha1big}[\tau^K_{\ell-1} \wedge T \wedge \theta_{k,m,C}^K,\tau^K_{\ell} \wedge T \wedge \theta_{k,m,C}^K]
=[\tau^K_{\ell-1},\tau^K_{\ell}].\end{equation}

Our aim now is to find the duration $\tau^K_{\ell}-\tau^K_{\ell-1}$ and to construct the set 
$M^K_{\ell}$ knowing the set $M^K_{\ell-1}$.

To reach $\tau^K_{\ell}$, two events are possible. Either one {living non resident} trait
reaches a size of order $K$, or a new mutant appears.

Let us consider the first type of event.
In fact, we have to be more precise on the time when a new trait has a size which reaches order $K$, this is 
why we defined $s_k^K$ {as} the time when one trait has a size which reaches order $K^{1-\eps_k}$. 
Notice that we may choose $\eps_k$ small enough to be sure that it corresponds to the trait {whose exponent reaches 1}
at time $s_k$ in the deterministic sequence $(s_j,j \in \N)$ defined in Theorem \ref{ThmConv}.
(if there exist two such traits, condition $(iv)(a)$ is fulfilled and $T_0$ is set to $s_k$).
Notice that if $f_{u,\v_{k-1}}<0$, for any $w \in V$,
\begin{equation} \label{decreasing} t \mapsto \beta^K_u(\tau^K_{\ell-1})+(t-\tau^K_{\ell-1})f_{u,\v_{k-1}}-\frac{d(u,w)}{\alpha}\end{equation}
is decreasing and thus will not reach $1-\eps_k$ if it is smaller than this value at time $\tau^K_{\ell-1}$.
Hence if we denote by $u_0$ the element of $M^K_{\ell-1}$ such that $\beta^K_{u_0}(\tau^K_{\ell})=1-\eps_k$, we get
$$ 1-\eps_k=\beta^K_{u_0}(\tau^K_{\ell})\cong\max_{\substack{u\in M^K_{\ell-1}
\\f_{u,\v_{k-1}}>0 }}\left[\beta^K_u(\tau^K_{\ell-1})+
(t-\tau^K_{\ell-1})f_{u,\v_{k-1}}-\frac{d(u,u_0)}{\alpha}\right].  $$
Now assume by contradiction that there is $u_1 \neq u_0 \in M^K_{\ell-1}$
such that:
\begin{align*} 1-\eps_k= \beta^K_{u_0}(\tau^K_{\ell})&\cong\beta^K_{u_1}(\tau^K_{\ell-1})+(\tau^K_{\ell}-\tau^K_{\ell-1})f_{u_1,\v_{k-1}}-\frac{d(u_1,u_0)}{\alpha}.
\end{align*}
This implies 
$$ \beta^K_{u_1}(\tau^K_{\ell})\geq \beta^K_{u_1}(\tau^K_{\ell-1})+(\tau^K_{\ell}-\tau^K_{\ell-1})f_{u_1,\v_{k-1}}>1, $$
as soon as $\eps_k<1/\alpha$,
which yields a contradiction. This implies that if there exists 
$u_0 \in M^K_{\ell-1}$ such that $\beta^K_{u_0}(\tau^K_{\ell})=1-\eps_k,$
then
\begin{align*}
 \beta^K_{u_0}(\tau^K_{\ell})&\cong\beta^K_{u_0}(\tau^K_{\ell-1})+(\tau^K_{\ell}-\tau^K_{\ell-1})f_{u_0,\v_{k-1}}
\end{align*}
and {with high probability}, the value of $\tau^K_\ell-\tau^K_{\ell-1}$ satisfies, 
\begin{equation} \label{first_case_taul} \tau^K_\ell-\tau^K_{\ell-1} \cong
\min_{\substack{w\in M^K_{\ell-1}:\\f_{w,\v_{k-1}}>0}}\frac{1-\beta^K_w(\tau^K_{\ell-1})}{f_{w,\v_{k-1}}}.\end{equation}

Let us now consider the second type of event, that is to say that there exist
$u_0 \notin M^K_{\ell-1}$ and $u_1 \in M^K_{\ell-1}$ such that $d(u_1,u_0)=1$ and $\beta^K_{u_1}(\tau^K_{\ell})=1/\alpha.$
Notice again than if $f_{u,\v_{k-1}}<0$, 
the function defined in \eqref{decreasing} is decreasing and thus will not reach $1/\alpha$ if it is smaller than this value at time $\tau^K_{\ell-1}$. 

By definition we have
\begin{align*} \frac{1}{\alpha}= \beta^K_{u_1}(\tau^K_{\ell})& \cong \max_{u\in M^K_{\ell-1}}\left[\beta^K_u(\tau^K_{\ell-1})+
(\tau^K_{\ell}-\tau^K_{\ell-1})f_{u,\v_{k-1}}-\frac{d(u,u_1)}{\alpha}\right].
\end{align*}
Denote by $u_2 \in M^K_{\ell-1}$ the trait realizing the maximum in the previous equation, that is to say
$$\frac{1}{\alpha}\cong \beta^K_{u_2}(\tau^K_{\ell-1})+(\tau^K_{\ell}-\tau^K_{\ell-1})f_{u_2,\v_{k-1}}-\frac{d(u_2,u_1)}{\alpha}.
$$

This equality can be rewritten {as}

$$\beta^K_{u_2}(\tau^K_{\ell-1})+(\tau^K_{\ell}-\tau^K_{\ell-1})f_{u_2,\v_{k-1}}\cong\frac{d(u_2,u_1)+1}{\alpha} .
$$
Let us now make a \textit{reductio ad absurdum} to prove that $d(u_2,u_1)+1=d(u_2,u_0)$. Let us thus assume that 
\begin{equation}\label{arg_abs}
 d(u_2,u_1)+1>d(u_2,u_0) \Leftrightarrow d(u_2,u_1)\geq d(u_2,u_0) ,
\end{equation}
and take $u_1'$ such that 
$$ d(u_2,u_1')+1=d(u_2,u_0). $$
Let us first assume (we will prove it later) that $u_1' \in M_{\ell-1}^K$. In this case, using the proof for the lower bound, we obtain that with high probability
\begin{align*} \beta^K_{u_1'}(\tau^K_{\ell}) &\geq \beta^K_{u_2}(\tau^K_{\ell-1})+(\tau^K_{\ell}-\tau^K_{\ell-1})f(u_2,\tilde{\v}_{k-1},k,-)-\frac{d(u_2,u_1')}{\alpha} \\
 & \geq  \beta^K_{u_2}(\tau^K_{\ell-1})+(\tau^K_{\ell}-\tau^K_{\ell-1})f(u_2,\tilde{\v}_{k-1},k,-)-\frac{d(u_2,u_1)-1}{\alpha}
 \cong \frac{2}{\alpha}.
\end{align*}
As $d(u_1',u_0)=1$, this means that $u_0$ becomes a living trait before the time $\tau_\ell^K$, which is in contradiction with the definition of $\tau_\ell^K$.

Let us now assume that $u_1' \notin M_{l-1}^K$ and consider a sequence of vertices $v_0=u_2,v_1,...,v_{d(u_2,u_1')}=u_1'$ such that 
$d(u_2,v_k)=k$ and $d(v_k,u_1')=d(u_2,u_1') -k $. Let 
$$ k_0:= \max\{0 \leq k \leq d(u_2,u_1')-1, v_k \in M^K_{\ell-1}\}. $$
Then 
$$ d(u_2,v_{k_0+1}) = d(u_2,u_1')-k_0-1 \leq d(u_2,u_1)-k_0-2, $$
and with high probability
\begin{align*} \beta^K_{v_{k_0+1}}(\tau^K_{\ell}) &\geq \beta^K_{u_2}(\tau^K_{\ell-1})+(\tau^K_{\ell}-
\tau^K_{\ell-1})f(u_2,\tilde{\v}_{k-1},k,-)-\frac{d(u_2,v_{k_0+1})}{\alpha} \\
 & \geq \beta^K_{u_2}(\tau^K_{\ell-1})+(\tau^K_{\ell}-\tau^K_{\ell-1})f(u_2,\tilde{\v}_{k-1},k,-)-\frac{d(u_2,u_1)-2}{\alpha} \cong  \frac{3}{\alpha},
\end{align*}
and thus $v_{k_0+1}$ becomes a living trait before the time $\tau_\ell^K$, which again is in contradiction with the definition 
of $\tau_\ell^K$.
We thus obtain a contradiction and deduce that \eqref{arg_abs} is not satisfied.
We conclude that
$$\beta^K_{u_2}(\tau^K_{\ell-1})+(\tau^K_{\ell}-\tau^K_{\ell-1})f_{u_2,\v_{k-1}}\cong\frac{d(u_2,u_1)+1}{\alpha}= \frac{d(u_2,u_0)}{\alpha}.
$$
Hence, when $\tau^K_{\ell}$ corresponds to the arrival of a new mutant,
\begin{align} \label{second_case_taul}
\tau^K_\ell-\tau^K_{\ell-1}
\cong \min_{\substack{w\in M^K_{\ell-1}
\\f_{w,\v_{\ell-1}}>0}}
\frac{\frac{d(w,V\backslash M^K_{\ell-1})}{\alpha}-\beta^K_w(\tau^K_{\ell-1})}{f_{w,\v_{k-1}}} . \notag\\
\end{align}

Combining \eqref{first_case_taul} and \eqref{second_case_taul}, we finally obtain:
\begin{align*}
\tau^K_\ell-\tau^K_{\ell-1}
&\cong\min_{\substack{w\in M^K_{\ell-1}:\\f_{w,\v_{k-1}}>0}}
\frac{\left(1\land\frac{d(w,V\backslash M^K_{\ell-1})}{\alpha}\right)-\beta^K_w(\tau^K_{\ell-1})}{f_{w,\v_{k-1}}}.
\end{align*}

To obtain $M^K_{\ell}$ from $M^K_{\ell-1}$, we suppress the traits $w \in M^K_{\ell-1}$ such that 
$ \beta^K_w(\tau^K_\ell)=0 $ (if condition $(iv)(c)$ is not satisfied, otherwise $T_0$ is set to $s_k$) 
and if $\tau_\ell \neq s_k$, we add the traits which are at distance 1 from the $w \in V$ satisfying
\begin{align*}
w\in\argmin_{\substack{w\in M^K_{\ell-1}:\\f_{w,\v_{k-1}}>0}}\frac{\left(1\land\frac{d(w,V\backslash M^K_{\ell-1})}
{\alpha}\right)-\beta^K_w(\tau^K_{\ell-1})}{f_{w,\v_{k-1}}}.
\end{align*}

\subsection{Value of $\theta_{k,m,C}^K$ and convergence of $s_k^K$ to $s_k$} \label{pageC}

Recall the definition of $\theta_{k,m,C}^K$ in \eqref{def_theta_k}.
We thus have constructed, on the time interval 
$[(\sigma_{k-1}^K \wedge T)\log K,(s_{k}^K \wedge \theta_{k,m,C}^K \wedge T)\log K]$, the times 
$(\tau_\ell^K, \ell \in \N)$ and the sets $(M^K_\ell, \ell \in \N)$ of living traits between times 
$\tau_\ell^K$ and $\tau_{\ell+1}^K$. We will now study the dynamics of the process on the time interval 
$[(\sigma_{k-1}^K \wedge T)\log K,(\sigma_{k}^K \wedge \theta_{k,m,C}^K\wedge T)\log K]$
($\sigma_k^K$ to be defined later in order to satisfy Assumption \ref{ass_sigma_k}).
Recall that $l_k^K$ is the trait $w \in V$ such that $\beta_w^K(s_k^K)=1-\eps_k$ and introduce 
$$\eta_k:= 2 \eps_k/\left(f_{l_k^K,\v_{k-1}}-\left(b_{l_k^K}+\left(\sum_{x \in V}c_{l_k^K,x}\right)(C+m)\right)\eps_k\right)$$

We will first prove that
\begin{align}\label{encadr_theta}
\lim_{K \to \infty} \P \left( s_k^K \leq \theta_{k,m,C}^K \leq s_k^K + \eta_k \Big| s_k^K < T \right)=1.
\end{align}
The first step consists in showing that 
\begin{equation}
 \lim_{K \to \infty} \P \left(\theta_{k,m,C}^K < s_k^K \Big| s_k^K < T \right)=0. \label{first_step}
\end{equation}
By definition of $s_k^K$, we have
$$ \sup_{w \in V \smallsetminus \v_{k-1}} \sup_{\sigma_{k-1}^K \leq t \leq s_k^K} \beta_w^K(t)\leq 1-\eps_k. $$
Moreover, applying Lemma \ref{lemmaC1} to $\v_{k-1}$ we obtain that 
$$ \lim_{K \to \infty} \P \left( \forall t \in [\sigma_{k-1}^K,s_k^K],
\sup_{w \in \v_{k-1}} \Big| \frac{N_w^K(t \log K)}{K} - \bar{n}_w(\v_{k-1}) \Big|\leq C \eps_k \Big| s_k^K < T \right)=1. $$
As a consequence, \eqref{first_step} holds true.
Notice that the value of $C$ in the definition of $\theta_{k,m,C}^K$ in \eqref{def_theta_k} is a consequence of the previous limit.
The constant $C$ is the one needed for Lemma \ref{lemmaC1} to hold, and thus depends on the parameters of the process.

Now assume by contradiction that
$$ s_k^K + \eta_k \leq \theta_{k,m,C}^K <T. $$
Then on the time interval $[s_k^K ,s_k^K + \eta_k]$, by definition of $\theta_{k,m,C}^K$, the $l^K_k$ population has a growth rate bounded from below by
$$f_{l_k^K,\v_{k-1}}-\left(b_{l_k^K}+\left(\sum_{x \in V}c_{l_k^K,x}\right)(C+m)\right)\eps_k.$$
Hence by coupling, with high probability, 
$$ \beta_{l_k^K}^K(s_k^K+\eta_k)\geq 1-\eps_k + 
\left(f_{l_k^K,\v_{k-1}}-\left(b_{l_k^K}+\left(\sum_{x \in V}c_{l_k^K,x}\right)(C+m)\right)\eps_k\right)\eta_k= 1+\eps_k, $$
which leads to a contradiction, as the total population size cannot be of order larger than $K$ in the limit $K\to\infty$, see Corollary \ref{betabound}.

This proves \eqref{encadr_theta}.
In particular, this implies that $s_k^K$ converges to $s_k$ in probability when $K$ goes to infinity, 
as soon as $T>s_k$.

\subsection{Value of the process at time $\theta_{k,m,C}^K \log K$}

We are now interested in the value of the process at time $\theta_{k,m,C}^K \log K$.
First notice that according to Proposition A.2 in \cite{CM11} and \eqref{encadr_theta},
\begin{multline} \label{stay_close}
\lim_{K \to \infty} \P \left( \forall t \in [\sigma_{k-1}^K,(s_k^K+\eta_k)\wedge \theta_{k,m,C}^K], w \in \v_{k-1},
\left|\frac{N_w^K(t \log K)}{K} - \bar{n}_w(\v_{k-1}) \right| <C\eps_k \right)\\
= \lim_{K \to \infty} \P \left( \forall t \in [\sigma_{k-1}^K,\theta_{k,m,C}^K], w \in \v_{k-1},\left|\frac{N_w^K(t \log K)}{K} - \bar{n}_w(\v_{k-1}) \right| <C\eps_k \right)=1. \end{multline}
Notice that $m$ has to be chosen small enough for this limit to hold, and thus depends on the parameters
of the Lotka-Volterra deterministic system associated to $\v_{k-1}$. 
To be more precise, $m$ has to be chosen small enough for the assumption \eqref{hyp_petites_dev} in Lemma \ref{lemmaC1} 
to hold true with $\eps_k$ in place of $\eps$.
We choose such an $m$ in the 
definition of $\theta_{k,m,C}^K$ in \eqref{def_theta_k}.

Hence we obtain that {with high probability},
$$ \sum_{w \in V \smallsetminus \v_{k-1}} N_w^K(\theta_{k,m,C}^K \log K) \geq m \eps_k K. $$

If condition $(iv)(a)$ of Theorem \ref{ThmConv} is satisfied $T_0$ is set at $s_k$ and the induction is stopped.
Otherwise there exists $\gamma>0$ such that if $\eps_k$ is small enough $\beta_w^K(s_k^K)<1-\gamma$ for every 
$w \in V \smallsetminus (\tilde{\v}^K_{k-1} \cup \{l^K_k\})$. Thus again by coupling, as the growth rates of the populations are 
limited 
and $\eta_k$ may be as small as we want, with high probability,
\begin{equation} \label{stay_small} 
\sum_{\substack{w \in V \smallsetminus \tilde{\v}^K_{k-1}, w \neq l^K_k }} 
N_w^K(\theta_{k,m,C}^K \log K) \leq K^{1-\gamma/2}.
 \end{equation}
From the two last inequalities we deduce that with high probability,
\begin{equation} \label{becomes_big}
N_{l^K_k}^K(\theta_{k,m,C}^K\log K) \geq m \eps_k K/2.
\end{equation}

\subsection{Construction of $\sigma_k^K$ and Assumption \ref{ass_sigma_k}}

Let us now introduce the stopping time $\sigma_k^K$, via:
\begin{equation} \label{def_sigma_k}
 \sigma_k^K := \inf \{ t \geq \theta_{k,m,C}^K, \forall w \in \v^K_{k}, |N^K_w(t \log K)/K-\bar{n}_w(\v_k)|\leq \eps_k\}.
\end{equation}

The last step of the proof consists in showing that $\sigma_k^K$ indeed satisfies Assumption \ref{ass_sigma_k}.
First $\sigma_k^K \log K $ is a stopping time. Second, from \eqref{stay_close}, \eqref{stay_small}, \eqref{becomes_big} and an application 
of Lemma \ref{lemmaC1}
there exists $T(\eps_k)<\infty$ such that
 $$ \lim_{K \to \infty} \P \left( \left| N_w^K(\theta_{k,m,C}^K \log K+T(\eps_k))/K-\bar{n}_w(\v_k) \right|
 \leq \eps_k, \forall w \in \v_k \right)=1. $$
Moreover, during a time of order one, the order of population sizes does not vary more than a constant times $\eps_k$
(result similar in spirit to Lemma B.9 in \cite{htrans19}). Adding that $s_k^K$ converges to $s_k$ in probability
when $K$ goes to infinity, as well as \eqref{encadr_theta}, we obtain that Assumption \ref{ass_sigma_k} holds. 
It ends the proof of Theorem \ref{ThmConv} and Proposition \ref{CorEquilibria}.


\appendix

\section{Couplings with branching processes and logistic processes with immigration}

The aim of this section is to collect various couplings of the populations with simpler processes like 
branching processes and logistic processes with immigration, and to state some properties of these simpler processes.
These results have been derived in \cite{htrans19} (note that we need to slightly generalise some of them), 
and we state them for the sake of readability.
For simplicity we keep the notations of \cite{htrans19}.

\subsection{Branching process} \label{def_BPK}

In this subsection, we recall Lemma A.1 of \cite{htrans19}, which describes the dynamics of a birth and death process on a $\log K$ time scale.
For $b,d,\beta \geq 0$, let $BP_K(b,d,\beta)$ denote the law of a process $(Z^K(t),t \geq 0)$ 
with initial state $Z^K(0)= \lfloor K^\beta -1 \rfloor$, individual birth rate $b$ and individual death rate $d$.

\begin{lemma}[Lemma A.1 in \cite{htrans19}] \label{lemmaA1}
 Let $(Z^K(t), t \geq 0)$ be a $BP_K (b, d, \beta)$ process such that $\beta>0$. The process
$(\log(1+Z^K(t \log K))/\log K ,t > 0)$
 converges when $K$ tends to infinity in probability in $L^\infty([0, T ])$ for all $T > 0$ to {the continuous deterministic function given by
 $$\bar{\beta} : t \mapsto \beta + (b-d)t \vee 0.$$
 }
 
 In addition, if $b<d$, for all $t> \beta / (d-b)$, 
 $$ \lim_{K \to \infty} \P \left( Z^K_{t \log K}=0 \right)=1. $$
 \end{lemma}

\subsection{Branching process with immigration} \label{sect_BPI}

In this subsection, we recall Lemma B.4 and Theorem B.5 of \cite{htrans19}, illustrated in Figure B.1 therein, which describe the dynamics of birth and death processes with immigration on a $\log K$ time scale.
For $b,d,\beta \geq 0$, $a,c \in \R$, $BPI_K(b,d,a,c,\beta)$ denotes the law of a process $(Z^K(t),t \geq 0)$ 
with initial state $Z^K(0)= \lfloor K^\beta -1 \rfloor$, individual birth rate $b$, individual death rate $d$, 
and immigration rate $K^ce^{as}$ at time $s \geq 0$.

\begin{lemma}[Lemma B.4 in \cite{htrans19}]\label{lemmaB4}
 Assume that $\beta<c$. Then for all $\eps>0$ and all $\bar{a}>|b-d| \vee |a| $, 
 $$ \lim_{K \to \infty} \P \left( Z^K(\eps \log K) \in \left[ K^{c-\bar{a}\eps},K^{c+\bar{a}\eps} \right] \right)=1. $$
\end{lemma}

\begin{lemma}[Theorem B.5 in \cite{htrans19}] \label{thmB5}
Let $(Z^K(t), t \geq 0)$ be a $BPI_K (b, d, a, c, \beta)$ process with $c \leq \beta$ and assume that
$\beta > 0$. The process
$(\log(1+Z^K(t \log K))/\log K ,t > 0)$
 converges when $K$ tends to infinity in probability in $L^\infty([0, T ])$ for all $T > 0$ to the continuous deterministic 
 function $\bar{\beta}$ given by
$$\bar{\beta} : t \mapsto (\beta + (b-d)t) \vee (c + at) \vee 0.$$
In addition, in the case where $c \neq 0$ or $a \neq 0$, for all compact intervals $I \subset \R_+$ which do
not intersect the support of $\bar{\beta}$,
$$\lim_{K \to \infty} \P \left(Z^K(t \log K) = 0, \forall t \in I \right)= 1.$$ 
\end{lemma}

	We will mostly use a corollary of those two lemmas, which is valid without the assumption $c\leq \beta$ but on a time interval 
	$[\delta,T]$, for any $\delta>0$. The idea of the proof has been explained in Remark \ref{rem-[eps,T]}.
	
	\begin{corollary} \label{thm-line-break}
	Let $(Z^K(t), t \geq 0)$ be a $BPI_K (b, d, a, c, \beta)$ process with $\beta\geq 0$, and either $c>0$ or both $c=0$ and $a>0$.
	For any $\delta>0$ and $T>0$, the process
	$(\log(1+Z^K(t \log K))/\log K ,t\in[\delta,T])$
	 converges when $K$ tends to infinity in probability in $L^\infty([\delta, T ])$ to the continuous deterministic 
	 function $\bar{\beta}$ given by
	$$\bar{\beta} : t \mapsto ((\beta\vee c) + (b-d)t) \vee (c + at) \vee 0.$$
\end{corollary}

\subsection{Logistic birth and death process with immigration}

We recall that for a subset $\v\subset V$ of traits that can coexist at a strictly positive equilibrium in the Lotka-Volterra
system \eqref{LV_system}, $\bar{n}(\v)\in\R_+^\v$ denotes this equilibrium. 
The next result states that if all traits in $\v$ have an initial population of order $K$ and 
the immigration of individuals with traits in $\v$ is small enough, the equilibrium $\bar{n}(\v) K$ is reached in 
a time of order $1$ and the populations of individuals whose traits belong to $\v$ will keep a size close to its 
equilibrium during a time of order larger than $\log K$

This result is a generalisation of Lemma C.1 in \cite{htrans19} to the multidimensional case and with (slightly) varying rates.

We thus consider a subset $\v\subset V$ of traits and denote by $(\mathbf{b}_{\v}(t), t \geq 0):=((b_w(t), w \in \v),t \geq 0)$, $(\mathbf{d}_{\v}(t), t \geq 0):=((d_w(t), w \in \v),t \geq 0)$, and $(\mathbf{c}_{\v}(t), t \geq 0):=((c_{w_1,w_2}(t) , (w_1,w_2) \in \v^2), t \geq 0)$ its birth, natural death, and death by competition rates that we allow to vary in time, as well 
as $(\mathbf{g}_{\v}(t), t \geq 0):=((g_w(t), w \in \v), t \geq 0 )$ a function with values in $\R_+^\v$.
We denote by $LBDI_K(\mathbf{b}_{\v},\mathbf{d}_{\v},\mathbf{c}_{\v},\mathbf{g}_{\v})$ the law of a logistic 
birth and death process with immigration $\mathbf{Z}^K:=((Z_w(t)^K, w \in \v), t \geq 0)$ where, at time $t$, an 
individual with a trait $w \in \v$ has a birth rate $b_w(t)$, a death 
rate $d_w(t) + \sum_{x \in \v} c_{w,x}(t) Z_x^K(t)/K $ and an immigration rate $g_w(t)$.

\begin{lemma}\label{lemmaC1}
Let $T>0$, $\v\subset V$ and assume that the mutation-free Lotka-Volterra system \eqref{LV_system} associated to $\v$ {and with rates
$(\bar{\mathbf{b}}_\v,\bar{\mathbf{d}}_\v,\bar{\mathbf{c}}_\v) \in (\R_*^+)^\v \times (\R_*^+)^\v \times (\R_*^+)^{\v^2}$} admits 
a unique positive globally attractive stable equilibrium $\bar{n}_w(\v)$.
Assume that $Z^K$ follows the law $LBDI_K(\mathbf{b}_{\v},\mathbf{d}_{\v},\mathbf{c}_{\v},\mathbf{g}_{\v})$ and 
that
\begin{equation} \label{hyp_petites_dev}
 \sup_{w_1,w_2\in \v} \left\{ \left|\mathbf{b}_{w_1}(t)-\bar{\mathbf{b}}_{w_1}\right|,\left|\mathbf{d}_{w_1}(t)-\bar{\mathbf{d}}_{w_1}\right|,
\left|\mathbf{c}_{w_1,w_2}(t)-\bar{\mathbf{c}}_{w_1,w_2}\right|  \right\}<\eps  
\end{equation}
and $g_w(t)\leq K^{1-\eta}$ for all $t \in [0, T \log K]$, $w \in \v$ for some 
$\eps,\eta>0$.
\begin{enumerate}
 \item[(i)] There exists $C,\eps_0>0$ such that if $\eps \leq \eps_0$ and $ \| \mathbf{Z}^K(0)/K- \bar{n}(\v) \|_\infty \leq \eps$, then
 $$ \lim_{K \to \infty} \P \left( \forall t \in [0,T\log K],  \| \mathbf{Z}^K(t)/K- \bar{n}(\v) \|_\infty \leq C\eps  \right)=1. $$
 \item[(ii)] For all $\eps_1,\eps_2>0$, there exists $T(\eps_1,\eps_2)<\infty$ such that for all initial condition 
 $\mathbf{Z}^K(0)$ such that $Z_v^K(0)/K \geq \eps_1 \ \forall v \in \v$ we have that
 $$ \lim_{K \to \infty} \P \left( \| \mathbf{Z}^K(T(\eps_1,\eps_2))/K- \bar{n}(\v) \|_\infty \leq \eps_2  \right)=1. $$
\end{enumerate}
\end{lemma}

{\begin{proof}
The case where the functions $\mathbf{b}_{\v}$, $\mathbf{d}_{\v}$, $\mathbf{c}_{\v}$ are constant is a direct generalisation of Lemma C.1 in \cite{htrans19}, whose proof
follows arguments similar to the ones given 
in \cite{Cha06,CM11} or in the Proposition 4.2 in \cite{champagnat2014adaptation} to handle the addition of (negligible) 
immigration. We do not provide it. Let us explain how we deal with varying rates for point (i). Let us choose $w_0 \in \v$, and introduce for $w_1,w_2 \in \v$:
$$ \tilde{b}_{w_1}= \left\{\begin{array}{ll}
                     \bar{b}_{w_1}+\eps & \text{if } w_1 \neq w_0\\
                     \bar{b}_{w_1}-\eps & \text{if } w_1 = w_0
                    \end{array} \right. $$
$$ \tilde{d}_{w_1}= \left\{\begin{array}{ll}
                     \bar{d}_{w_1}-\eps & \text{if } w_1 \neq w_0\\
                     \bar{d}_{w_1}+\eps & \text{if } w_1 = w_0
                    \end{array} \right. $$
$$ \tilde{c}_{w_1,w_2}= \left\{\begin{array}{ll}
                     \bar{c}_{w_1,w_2}-\eps & \text{if } w_1 \neq w_0\\
                     \bar{c}_{w_1,w_2}+\eps & \text{if } w_1 = w_0
                    \end{array} \right. $$
Then we can couple a process $Z^K$ with the law $LBDI_K(\bar{\mathbf{b}}_{\v},\bar{\mathbf{d}}_{\v},\bar{\mathbf{c}}_{\v},\mathbf{g}_{\v})$ with a process 
$\tilde{Z}^K$ with the law $LBDI_K(\tilde{\mathbf{b}}_{\v},\tilde{\mathbf{d}}_{\v},\tilde{\mathbf{c}}_{\v},\mathbf{g}_{\v})$ such that for every $t \geq 0$, 
$\tilde{Z}_{w_0}^K(t) \leq Z^K_{w_0}(t) $ and $\tilde{Z}_{w}^K(t) \geq Z^K_{w}(t) $ for every $w \in \v \setminus w_0$. Moreover, as the equilibrium of a Lotka-Volterra system is continuous 
with respect to its coefficient, there is a positive $\tilde{C}$ such that for $\eps$ small enough, and if we denote by $\bar{n}^{(w_0)}(\v)$ the equilibrium of the Lotka-Volterra system with 
the coefficients $\tilde{\mathbf{b}}_{\v},\tilde{\mathbf{d}}_{\v},\tilde{\mathbf{c}}_{\v}$ we have just introduced, $\|\bar{n}^{(w_0)}(\v)-\bar{n}^{(w_0)}(\v)\| \leq \tilde{C} \eps$.
Hence applying the point (i) for the process $\tilde{Z}^K$, we obtain upper bounds for coordinates $w \neq w_0$ and a lower bound for the coordinate $w_0$, for the process $Z^K$. 
Doing the same and the reverse bounds for the other elements of $\v$ gives the result for some $C>\tilde{C}$ that takes into account the fluctuations around the varied equilibria.
\end{proof}}

We end this section with a result stating that the time needed for the total population size of a logistic birth and death process (with or without 
mutations) to reach (and stay smaller than) an order $K$ is of order one for $K$ large enough.

\begin{corollary}\label{betabound}
Let us consider a subset $\v\subset V$ of traits, $(\mathbf{b}_\v,\mathbf{d}_\v,\mathbf{c}_\v)$ be in $(\R_*^+)^\v \times (\R_*^+)^\v \times (\R_*^+)^{\v^2}$ and 
let $Z^K$ follow the law $LBDI_K(\mathbf{b}_{\v},\mathbf{d}_{\v},\mathbf{c}_{\v},0)$, and $\mathcal{Z}^K$ denote the total 
population size of the process $Z^K$. For every $\eps>0$ there exists $T(\eps)<\infty)$ such that for $t>T(\eps)$
\begin{align*}
\lim_{K\to\infty}\P\left(\frac{\log(1+\mathcal{Z}^K(t))}{\log K}<1+\eps\right)=1.
\end{align*}
\end{corollary}

\begin{remark}
Notice that this result only treats mutation-free logistic birth and death processes. However, mutation within $\v$ does not affect the total population size and hence the result can be transferred to such cases. Considering $\v=V$, Corollary \ref{betabound} therefore implies the same asymptotic bound for the total population size of the process that we consider in Theorem \ref{ThmConv} and Proposition \ref{CorEquilibria}, and hence also for each subpopulation.
\end{remark}

\begin{proof}
The process $\mathcal{Z}^K$ increases by $1$ at a rate
$$ \sum_{w \in \v} b_w Z_w^K \leq (\sup_{w \in \v}b_w)\mathcal{Z}^K=:\mathcal{B}\mathcal{Z}^K $$
and decreases by $1$ at a rate
\begin{align*} \sum_{w \in \v} \left(d_w + \sum_{u \in \v}\frac{c_{w,u}}{K}Z_u^K \right) Z_w^K & \geq \frac{1}{K} (\inf_{u \in \v}c_{u,u}) \sum_{w \in \v}\left(Z_w^K\right)^2 \\
 &  \geq \frac{1}{K}(\inf_{u \in \v}c_{u,u}) \frac{1}{Card ( \v)}\left(\mathcal{Z}^K\right)^2 =:\frac{\mathcal{C}}{K} \left(\mathcal{Z}^K\right)^2.
\end{align*}
Hence the process $\mathcal{Z}^K$ can be coupled with a logistic birth and death process $\mathcal{N}^K$ with individual birth rate $\mathcal{B}$
and individual death rate $\mathcal{C} \mathcal{N}^K/K$ in such a way that for every $t \geq 0$, if $ \mathcal{Z}^K(0)= \mathcal{N}^K(0)$
$$ \mathcal{Z}^K(t) \leq \mathcal{N}^K(t). $$
But from Chapter 11, Theorem 2.1 in \cite{ethiermarkov}, we know that on any finite time interval, the rescaled process $\mathcal{N}^K/K$ converges in probability
to the solution to the logistic equation $\dot{\varkappa}= \varkappa(\mathcal{B}-\mathcal{C}\varkappa)$, $\varkappa(0)=\varkappa_0$ if 
$\mathcal{N}^K(0)/K$ converges in probability to $\varkappa_0$. The one dimensional logistic equation has an explicit solution, and in 
particular, we know that its equilibrium is $\mathcal{B}/\mathcal{C}$, that it comes down from infinity, and that it takes a time 
$$\frac{1}{\mathcal{B}} \log \left( \frac{\bar{\varkappa}}{\bar{\varkappa}-\mathcal{B}/\mathcal{C} } \right) $$
to reach $\bar{\varkappa}>\mathcal{B}/\mathcal{C}$ from an infinite initial condition. As a consequence, $\mathcal{N}^K$ takes a time of order one 
to become smaller than $2 \bar{\varkappa} K$, and as $\mathcal{B}/\mathcal{C}$ is a globally hyperbolic equilibrium for the function $\varkappa$, classical 
large deviation results (see \cite{dupuis2011weak} for instance) entail that $\mathcal{N}^K/K$ will stay an exponential (in $K$) time in any compact interval of $R_+^*$
including $\mathcal{B}/\mathcal{C}$. This concludes the proof.
\end{proof}

\bibliographystyle{abbrv}
\bibliography{Library-Biomaths}

\begin{thebibliography}{10}

\bibitem{BerBruShi2016}
J.~Berestycki, E.~Brunet, and Z.~Shi.
\newblock {The number of accessible paths in the hypercube}.
\newblock {\em Bernoulli}, 22(2):653--680, may 2016.

\bibitem{BilSma15}
S.~Billiard and C.~Smadi.
\newblock {The interplay of two mutations in a population of varying size: a
  stochastic eco-evolutionary model for clonal interference}.
\newblock {\em Stoch. Process. their Appl.}, 127(3):701--748, 2017.

\bibitem{BilSma2019}
S.~Billiard and C.~Smadi.
\newblock {Stochastic dynamics of three competing clones: Conditions and times
  for invasion, coexistence, and fixation}.
\newblock {\em Am. Nat.}, oct 2019.

\bibitem{BolPac1}
B.~Bolker and S.~W. Pacala.
\newblock {Using Moment Equations to Understand Stochastically Driven Spatial
  Pattern Formation in Ecological Systems}.
\newblock {\em Theor. Popul. Biol.}, 52(3):179--197, 1997.

\bibitem{BolPac2}
B.~M. Bolker and S.~W. Pacala.
\newblock {Spatial moment equations for plant competition: understanding
  spatial strategies and the advantages of short dispersal}.
\newblock {\em Am Nat}, 153(6):575--602, 1999.

\bibitem{BovCoqSma2018}
A.~Bovier, L.~Coquille, and C.~Smadi.
\newblock {Crossing a fitness valley as a metastable transition in a stochastic
  population model}.
\newblock {\em Ann. Appl. Probab.}, 29(6):3541--3589, dec 2019.

\bibitem{Cha06}
N.~Champagnat.
\newblock {A microscopic interpretation for adaptive dynamics trait
  substitution sequence models}.
\newblock {\em Stoch. Process. their Appl.}, 116(8):1127--1160, 2006.

\bibitem{champagnat2014adaptation}
N.~Champagnat, P.-E. Jabin, and S.~M{\'{e}}l{\'{e}}ard.
\newblock {Adaptation in a stochastic multi-resources chemostat model}.
\newblock {\em J. Math. Pures Appl.}, 101(6):755--788, jun 2014.

\bibitem{CM11}
N.~Champagnat and S.~M{\'{e}}l{\'{e}}ard.
\newblock {Polymorphic evolution sequence and evolutionary branching}.
\newblock {\em Probab. Theory Relat. Fields}, 151(1-2):45--94, oct 2011.

\bibitem{htrans19}
N.~Champagnat, S.~M{\'{e}}l{\'{e}}ard, and V.~C. Tran.
\newblock {Stochastic analysis of emergence of evolutionary cyclic behavior in
  population dynamics with transfer}.
\newblock {\em ArXiv Prepr. 1901.02385}, 2019.

\bibitem{CMM13}
P.~Collet, S.~M{\'{e}}l{\'{e}}ard, and J.~A.~J. Metz.
\newblock {A rigorous model study of the adaptive dynamics of Mendelian
  diploids}.
\newblock {\em J. Math. Biol.}, 67(3):569--607, 2013.

\bibitem{Coron2017}
C.~Coron, M.~Costa, H.~Leman, and C.~Smadi.
\newblock {A stochastic model for speciation by mating preferences}.
\newblock {\em J. Math. Biol.}, 76(6):1421--1463, may 2018.

\bibitem{Cow2006}
M.~C. Cowperthwaite, J.~J. Bull, and L.~A. Meyers.
\newblock {From Bad to Good: Fitness Reversals and the Ascent of Deleterious
  Mutations}.
\newblock {\em PLoS Comput. Biol.}, 2(10):e141, 2006.

\bibitem{de2014empirical}
J.~A.~G. De~Visser and J.~Krug.
\newblock Empirical fitness landscapes and the predictability of evolution.
\newblock {\em Nature Reviews Genetics}, 15(7):480--490, 2014.

\bibitem{DePristo2007}
M.~A. DePristo, D.~L. Hartl, and D.~M. Weinreich.
\newblock {Mutational Reversions During Adaptive Protein Evolution}.
\newblock {\em Mol. Biol. Evol.}, 24(8):1608--1610, aug 2007.

\bibitem{DieLaw}
Dieckmann and Law.
\newblock {Moment approximations of individual-based models}.
\newblock In U.~Dieckmann, R.~Law, and J.~A.~J. Metz, editors, {\em Geom. Ecol.
  Interact. Simpl. Spat. Complex.}, pages 252--270. Cambridge University Press,
  2000.

\bibitem{DL96}
U.~Dieckmann and R.~Law.
\newblock {The dynamical theory of coevolution: a derivation from stochastic
  ecological processes}.
\newblock {\em J. Math. Biol.}, 34(5-6):579--612, may 1996.

\bibitem{dupuis2011weak}
P.~Dupuis and R.~S. Ellis.
\newblock {\em {A Weak Convergence Approach to the Theory of Large
  Deviations}}.
\newblock Wiley Series in Probability and Statistics. John Wiley \& Sons, Inc.,
  Hoboken, NJ, USA, feb 1997.

\bibitem{DurMay2011}
R.~Durrett and J.~Mayberry.
\newblock {Traveling waves of selective sweeps}.
\newblock {\em Ann. Appl. Probab.}, 21(2):699--744, apr 2011.

\bibitem{ethiermarkov}
S.~N. Ethier and T.~G. Kurtz.
\newblock {Markov Processes: Characterization and Convergence}, 1986.

\bibitem{FM04}
N.~Fournier and S.~M{\'{e}}l{\'{e}}ard.
\newblock {A microscopic probabilistic description of a locally regulated
  population and macroscopic approximations}.
\newblock {\em Ann. Appl. Probab.}, 14(4):1880--1919, 2004.

\bibitem{GMK1997}
S.~A.~H. Geritz, J.~A.~J. Metz, {\'{E}}.~Kisdi, and G.~Mesz{\'{e}}na.
\newblock {Dynamics of Adaptation and Evolutionary Branching}.
\newblock {\em Phys. Rev. Lett.}, 78(10):2024--2027, mar 1997.

\bibitem{KL1987}
S.~Kauffman and S.~Levin.
\newblock {Towards a general theory of adaptive walks on rugged landscapes}.
\newblock {\em J. Theor. Biol.}, 128(1):11--45, sep 1987.

\bibitem{BovKra2018}
A.~Kraut and A.~Bovier.
\newblock {From adaptive dynamics to adaptive walks}.
\newblock {\em J. Math. Biol.}, 79(5):1699--1747, oct 2019.

\bibitem{Krug-rev-access-2019}
J.~Krug.
\newblock {Accessibility percolation in random fitness landscapes}.
\newblock {\em to Appear "Probabilistic Struct. Evol. ed. by E. Baake A.
  Wakolbinger}, mar 2019.

\bibitem{Len2003}
R.~E. Lenski, C.~Ofria, R.~T. Pennock, and C.~Adami.
\newblock {The evolutionary origin of complex features}.
\newblock {\em Nature}, 423(6936):139--144, may 2003.

\bibitem{maddamsetti2015adaptation}
R.~Maddamsetti, R.~E. Lenski, and J.~E. Barrick.
\newblock {Adaptation, clonal interference, and frequency-dependent
  interactions in a long-term evolution experiment with escherichia coli}.
\newblock {\em Genetics}, 200(2):619--631, jun 2015.

\bibitem{MS1970}
J.~{Maynard Smith}.
\newblock {Natural selection and the concept of a protein space}.
\newblock {\em Nature}, 225(5232):563--564, 1970.

\bibitem{MG96}
J.~A.~J. Metz, S.~A.~H. Geritz, G.~Meszena, F.~J.~A. Jacobs, and J.~S. van
  Heerwaarden.
\newblock {Adaptive Dynamics: A Geometrical Study of the Consequences of Nearly
  Faithful Reproduction}.
\newblock Iiasa working paper, IIASA, Laxenburg, Austria, 1995.

\bibitem{Motter2010}
A.~E. Motter.
\newblock {Improved network performance via antagonism: From synthetic rescues
  to multi-drug combinations}.
\newblock {\em BioEssays}, 32(3):236--245, mar 2010.

\bibitem{NK2011}
J.~Neidhart and J.~Krug.
\newblock {Adaptive walks and extreme value theory}.
\newblock {\em Phys. Rev. Lett.}, 107(17), oct 2011.

\bibitem{BovNeu}
R.~Neukirch and A.~Bovier.
\newblock {Survival of a recessive allele in a Mendelian diploid model}.
\newblock {\em J. Math. Biol.}, pages 1--54, nov 2016.

\bibitem{NovKru2015}
S.~Nowak and J.~Krug.
\newblock {Analysis of adaptive walks on NK fitness landscapes with different
  interaction schemes}.
\newblock {\em J. Stat. Mech. Theory Exp.}, 2015(6), jun 2015.

\bibitem{Orr2003}
H.~A. Orr.
\newblock {A minimum on the mean number of steps taken in adaptive walks}.
\newblock {\em J. Theor. Biol.}, 220(2):241--247, 2003.

\bibitem{PloKud}
J.~B. Plotkin and G.~Kudla.
\newblock {Synonymous but not the same: The causes and consequences of codon
  bias}, jan 2011.

\bibitem{SahMot2011}
S.~Sahasrabudhe and A.~E. Motter.
\newblock {Rescuing ecosystems from extinction cascades through compensatory
  perturbations}.
\newblock {\em Nat. Commun.}, 2(1), 2011.

\bibitem{SchKru2014}
B.~Schmiegelt and J.~Krug.
\newblock {Evolutionary Accessibility of Modular Fitness Landscapes}.
\newblock {\em J. Stat. Phys.}, 154(1-2):334--355, 2014.

\bibitem{smadi2015eco}
C.~Smadi.
\newblock {An eco-evolutionary approach of adaptation and recombination in a
  large population of varying size}.
\newblock {\em Stoch. Process. their Appl.}, 125(5):2054--2095, may 2015.

\bibitem{Smadi2017}
C.~Smadi.
\newblock {The Effect of Recurrent Mutations on Genetic Diversity in a Large
  Population of Varying Size}.
\newblock {\em Acta Appl. Math.}, 149(1):11--51, 2017.

\bibitem{Krug-rev-emp-2013}
I.~G. Szendro, M.~F. Schenk, J.~Franke, J.~Krug, and J.~A.~G. {De Visser}.
\newblock {Quantitative analyses of empirical fitness landscapes}.
\newblock {\em J. Stat. Mech. Theory Exp.}, 2013(1), jan 2013.

\bibitem{Zee1993}
M.~L. Zeeman.
\newblock {Hopf bifurcations in competitive three-dimensional
  {L}otka-{V}olterra systems}.
\newblock {\em Dynam. Stab. Syst.}, 8(3):189--217, 1993.

\end{thebibliography}

\end{document}